\documentclass[11pt]{amsart}
\usepackage{amsmath,amsthm, amscd, amssymb, amsfonts}
\usepackage[all]{xy}

\newcommand{\morf}{{\mathbf \psi}}

\DeclareMathOperator*{\prode}{\boxtimes}
\DeclareMathOperator*{\botimes}{\mbox{\boldmath $\otimes$}}
\newcommand{\sou}{\mathfrak s}
\newcommand{\tgt}{\mathfrak e}

\newcommand{\equal}[1]{
{\stackrel{{\textstyle #1}}{\ {\textstyle =}\ } }}
\newcommand{\onto}{\twoheadrightarrow} 
\newcommand{\map}[1]{\xrightarrow{#1}}

\newcommand{\abs}[1]{\lvert#1\rvert}   
\newcommand{\Abs}[1]{\bigl\vert#1\bigr\vert}
\newcommand{\norm}[1]{\lVert#1\rVert}
\newcommand{\Norm}[1]{\bigl\Vert#1\bigr\Vert}

\newcommand{\R}{{\mathcal R}}
\newcommand{\Rs}{\overline{\mathcal R}}

\newcommand{\frake}{{\mathfrak e}}
\newcommand{\fraks}{{\mathfrak s}}

\newcommand{\ffl}{f\!\ell}
\newcommand{\Rep}{\operatorname{Rep}}
\newcommand{\Repk}{\Rep_{\ku}}
\newcommand{\Quiv}{\operatorname{Quiv}}
\newcommand{\Bimod}{\operatorname{Bimod}}

\newcommand{\ku}{\Bbbk}
\newcommand{\K}{{\mathcal K}}

\newcommand{\N}{{\mathcal N}}
\newcommand{\G}{{\mathcal G}}

\newcommand{\pr}{\operatorname{pr}}

\newcommand{\Z}{{\mathbb Z}}

\newcommand{\F}{{\mathcal F}}
\newcommand{\C}{{\mathcal C}}
\newcommand{\D}{{\mathcal D}}
\newcommand{\Ec}{{\mathcal E}}
\newcommand{\uno}{{\bf 1}}

\newcommand{\Hc}{{\mathcal H}}
\newcommand{\Vc}{{\mathcal V}}
\newcommand{\Pc}{{\mathcal P}}

\newcommand{\Hb}{{\mathbb H}}
\newcommand{\Vb}{{\mathbb V}}

\newcommand{\Bg}{{\mathfrak B}}

\newcommand{\Hg}{{\mathfrak H}}
\newcommand{\Vg}{{\mathfrak V}}

\newcommand{\Ss}{{\mathcal S}}

\newcommand{\End}{\operatorname{End}}
\newcommand{\Lin}{\operatorname{Lin}}

\newcommand{\Res}{\operatorname{Res}}

\newcommand{\fpc}{\Quiv(\Pc)}

\newcommand{\fde}{\,{\triangleright}\,}
\newcommand{\fiz}{\,{\triangleleft}\,}

\DeclareMathOperator*{\Tim}{\times}
\newcommand{\Times}[2]{\sideset{_#1}{_#2}\Tim}

\numberwithin{equation}{section}\theoremstyle{plain}

\newtheorem{theorem}{Theorem}[section]
\newtheorem{lema}[theorem]{Lemma}
\newtheorem{lemma}[theorem]{Lemma}

\newtheorem{corollary}[theorem]{Corollary}

\newtheorem{prop}[theorem]{Proposition}
\newtheorem{proposition}[theorem]{Proposition}

\theoremstyle{definition}
\newtheorem{definition}[theorem]{Definition}
\newtheorem{exa}[theorem]{Example}

\newtheorem{claim}{Claim}

\newtheorem{obs}[theorem]{Remark}

\newcommand\id{\operatorname{id}}

\def\pf{\begin{proof}}
\def\epf{\end{proof}}

\theoremstyle{remark}

\begin{document}

\renewcommand{\baselinestretch}{1.2}
\thispagestyle{empty}
\title[Representations of matched pairs of groupoids]{Representations of matched pairs of groupoids
and applications to weak Hopf algebras}
\author{Marcelo Aguiar and Nicol\'as Andruskiewitsch}
\address{Department of Mathematics \newline
\indent Texas A\&M University \newline \indent
College Station, TX 77843-3368}
\email{maguiar@math.tamu.edu, \quad \emph{URL:}\/
http://www.math.tamu.edu/$\sim${}maguiar}
\address{Facultad de Matem\'atica, Astronom\'\i a y F\'\i sica
\newline \indent
Universidad Nacional de C\'ordoba
\newline
\indent CIEM -- CONICET
\newline
\indent (5000) Ciudad Universitaria, C\'ordoba, Argentina}
\email{andrus@mate.uncor.edu, \quad \emph{URL:}\/
http://www.mate.uncor.edu/andrus}
\thanks{First author supported in part by NSF grant DMS-0302423}
\thanks{Second author  supported in part by CONICET, Agencia C\'ordoba Ciencia, ANPCyT, and Secyt (UNC)}
\thanks{This paper is in final form and no version of it will be submitted elsewhere}
\date{February 7, 2004}
\begin{abstract} We introduce the category of set-theoretic representations of a matched
pair of groupoids. This is a monoidal category endowed with a monoidal functor $\ffl$ to the category
of quivers over the common base of the groupoids in the matched pair.
We study monoidal functors between two such categories of representations which preserve the functor
$\ffl$. We show that the centralizer of such a monoidal functor is the category of representations of a new
matched pair, which we construct explicitly. We introduce the notions of {\em double} 
of a matched pair of groupoids and \emph{generalized double} of a morphism of matched pairs. 
We show that the centralizer of
$\ffl$ is the category of representations of the dual matched pair, and
the centralizer of the identity functor (the center) is the category of representations of the double.
We use these constructions to classify the
braidings in the category of representations of a matched pair. Such braidings are parametrized by
certain groupoid-theoretic structures which we call {\em matched pairs of rotations}.
 Finally, we express our results in terms of the weak Hopf
algebra associated to a matched pair of groupoids. A matched pairs of rotations gives
rise to a quasitriangular structure for the associated weak Hopf algebra. The Drinfeld double of the weak
Hopf algebra of a matched pair is the weak Hopf algebra associated to the double matched pair. 
\end{abstract}

\maketitle

\section*{Introduction} 

An exact factorization of a group $\Sigma$ is a pair of subgroups $F$, $G$, such that
$\Sigma=FG$ and $F\cap G=\{1\}$. Such a factorization may be described without reference to the group
$\Sigma$, by means of a pair of actions of the factors $F$ and $G$ on each other satisfying certain
axioms: a matched pair of groups~\cite{k,mackey,tak0}.

Exact factorizations $\D = \Vc \Hc$ of groupoids are equivalent to structures of matched pairs
of groupoids on the pair $\Vc$, $\Hc$, and to structures of vacant double groupoids
~\cite{mk1,AN}. In this paper we work with matched pairs of groupoids, while
the notions from double category theory are implicit in our notations and illustrations.
Thus, a matched pair of groupoids $(\Vc,\Hc)$ consists of a pair of groupoids $\Vc$ and $\Hc$
with a common base (set of objects) $\Pc$
together with actions $\fde$ and $\fiz$ of each of them on the other satisfying certain simple axioms
(Definition~\ref{matchpairgpds}). This is summarized by the following illustration
\[\xymatrix@C+5pt@R-5pt{
{P}\ar@/_/[d]_{x\fde g}\ar@/^/[r]^{x} & {Q}\ar@/^/[d]^{g}\\  
{R}\ar@/_/[r]_{x\fiz g} & {S} }\]
where $g$ is an element (arrow) of the {\em vertical} groupoid $\Vc$ and $x$ is an element of the {\em
horizontal} groupoid $\Hc$.

A quiver over $\Pc$ is a set $\Ec$ equipped with two maps $p,q:\Ec\to\Pc$.
We introduce the notion of representations of a matched pair of groupoids. A representation
of $(\Vc,\Hc)$ is a quiver $\Ec$ over $\Pc$ together with an action of $\Hc$ and a grading over $\Vc$ which
are compatible (Definition~\ref{D:representation}). The matched pair structure allows us to construct
a monoidal structure on the category $\Rep(\Vc,\Hc)$ of representations. The forgetful functor
$\ffl:\Rep(\Vc,\Hc)\to\Quiv(\Pc)$ to the category of quivers over $\Pc$ is monoidal.

We introduce a notion of morphisms between matched pairs 
(Definition~\ref{D:morphism}). Associated to a morphism
$(\alpha,\beta):(\Vb,\Hb)\to(\Vc,\Hc)$ there is a monoidal functor
$\Res^\beta_\alpha:\Rep(\Vc,\Hc)\to\Rep(\Vb,\Hb)$, called the restriction
along $(\alpha,\beta)$.
Any restriction functor
preserves the forgetful functors $\ffl$ (Section~\ref{S:restriction}).
Our first main result states that any monoidal functor $\Rep(\Vc,\Hc)\to\Rep(\Vb,\Hb)$ which
preserves $\ffl$ must be a restriction functor (Theorem~\ref{restr-uniq}).

In Section~\ref{S:centers} we undertake the study of centralizers of such monoidal functors.
This includes the center of a category $\Rep(\Vc,\Hc)$ as the particular case when the
monoidal functor is the identity. Our second main result (Theorem~\ref{T:centralizer})
states that any such centralizer is again the category of representations of a (new) matched pair,
which we call a generalized double.  In particular, the center of $\Rep(\Vc,\Hc)$ is the category of
representations of the double $D(\Vc,\Hc)$ of $(\Vc,\Hc)$, and the centralizer of
$\ffl:\Rep(\Vc,\Hc)\to\Quiv(\Pc)$ is the category of representations of the dual of $(\Vc,\Hc)$
(Corollaries~\ref{C:center} and~\ref{C:dual}).

Doubles and generalized doubles are introduced and studied
in Section~\ref{S:morphisms}. The explicit description is given in
Theorem~\ref{T:double} and the functoriality of the construction in
Proposition~\ref{P:functoriality}. We also show that the constructions
of generalized doubles and that of duals commute 
(Proposition~\ref{P:gendouble-dual}).

The classification of braidings in the monoidal category $\Rep(\Vc,\Hc)$ is accomplished
in Section~\ref{S:braidings}. By general results (recalled in Section~\ref{S:generalities}),
braidings on $\Rep(\Vc,\Hc)$ are in bijective correspondence with monoidal sections of the canonical
functor $Z\bigl(\Rep(\Vc,\Hc)\bigr)\to\Rep(\Vc,\Hc)$. The latter is the restriction functor
along a canonical morphism of matched pairs $(\iota,\pi):(\Vc,\Hc)\to D(\Vc,\Hc)$.
Theorem~\ref{restr-uniq} reduces the
classification to the description of all sections of $(\iota,\pi)$. 
The result is that braidings are in bijective correspondence with {\em matched pairs of rotations}
for $(\Vc,\Hc)$. 
A rotation for $(\Vc,\Hc)$ is a morphism $\eta$ from the vertical groupoid $\Vc$ to the
horizontal groupoid $\Hc$ which is compatible with the actions $\fde$ and $\fiz$
(Definition~\ref{D:rotation}). Suppose $\eta$ and $\xi$ are rotations. Starting from a 
 pair of composable arrows $(g,f)$ in $\Vc$, one may rotate $g$ with $\eta$ and $f^{-1}$ with $\xi$, and
build the following diagram:
\[\xymatrix@R-5pt@C+15pt{
& {P}\ar@/_/[d]_{g} & {R}\ar@/_/[l]_{\xi(f)^{-1}\fiz g^{-1}} \\
{P}\ar@/_/[r]^{\eta(g)}\ar@/_/[d]_{\eta(g)\fde f} &{Q}\ar@/_/[d]_{f}\ar@/_/[u]_{g^{-1}} &
{S}\ar@/_/[l]^{\xi(f)^{-1}}\ar@/_/[u]_{\xi(f)^{-1}\fde g^{-1}}\\ 
{R}\ar@/_/[r]_{\eta(g)\fiz f} & {S}\ar@/_/[u]_{f^{-1}} & \\
}\]
There are then two ways of going from $P$ to $S$ via vertical arrows (through $Q$ or
through $R$). The pair $(\eta,\xi)$ is a matched pair of
rotations if these two ways coincide (Definition~\ref{D:R-matrix}). The correspondence between
braidings and matched pairs of rotations is settled in our third main result,
Theorem~\ref{T:classification}.

{}From a matched pair of finite groups one may construct a Hopf algebra $\ku(F,G)$, Takeuchi's
{\em bismash product}~\cite{k, tak0, maj}, which is part of a Hopf algebra extension
\[1 \to \ku^F \to  \ku(F,G) \to \ku G \to 1\,.\]

On the other hand, a construction of a weak Hopf algebra (or quantum groupoid)
$\ku(\Vc,\Hc)$ out of a matched pair $(\Vc,\Hc)$ of finite groupoids appears in~\cite{AN}. 
This fits in an extension
\[1 \to \ku^{\Vc} \to  \ku(\Vc,\Hc) \to \ku {\Hc} \to 1\,.\]
We review this construction in Section~\ref{S:weak}. There is a monoidal functor from the category
of representations of a matched pair $(\Vc,\Hc)$ to the category of modules over the weak Hopf algebra
$\ku(\Vc,\Hc)$  (the linearization functor).
In Theorem~\ref{T:positiveR}, we show that a matched pair of rotations $(\eta,\xi)$ for $(\Vc,\Hc)$ gives
rise to an quasitriangular structure $\R$ for $\ku(\Vc,\Hc)$, in the sense of~\cite{nik-v}. Explicitly,
\[\R=\sum \bigl(\xi(f)^{-1}\fiz g^{-1},g\bigr)
\otimes \bigl(\eta(g),f\bigr)\,,\]
the sum being over the pairs $(f,g)$ of composable arrows in the groupoid $\Vc$. We also provide an
explicit formula for the corresponding Drinfeld element.

The construction of a weak Hopf algebra from a matched pair commutes with
duals~\cite[Prop. 3.11]{AN}. We review this fact in Section~\ref{S:double-weak}, and
we show that the construction also commutes with doubles. 
Theorem~\ref{T:drinfeld} provides an explicit description for the
Drinfeld double of the weak Hopf algebra of a matched pair. In the case of matched pairs of groups,
the Drinfeld double was calculated in~\cite{mbg}.

The quasitriangular structures for $\ku(\Vc,\Hc)$ that we find specialize for the case of groups
to those constructed by Lu, Yan, and Zhu in~\cite{lyz2}. Our matched pairs of rotations
are in this case the {\em LYZ pairs} of~\cite{tak}. This reference introduces a notation reminiscent of
double categories, similar to ours. A different construction of these quasitriangular structures, still
for the case of groups, was given by Masuoka (see~\cite{tak}). 
Our approach, based on the calculation of the center of the category of representations, is
different from both~\cite{lyz2} and~\cite{tak}. 

In~\cite{lyz1,lyz2} it is shown that every Hopf algebra with a {\em positive} basis
is of the from $\ku(F,G)$ for some matched pair of groups $(F,G)$, and that all
positive quasitriangular structures for $\ku(F,G)$ are of the form mentioned above. We leave for future
work the question of whether these remain true for the case of groupoids.

The results of this paper are used in~\cite{A}.

\subsection*{Acknowledgements} We thank Jos\'e Antonio de la Pe\~na and the other organizers
of the {\em XV Coloquio Latinoamericano de \'Algebra}, held in Cocoyoc, M\'exico in July, 2004,
where our collaboration began.

\section{Matched pairs of groupoids}

\subsection{Groupoids}\label{S:basic}

\

Let $\Pc$ be a set. Given maps
$p:X\to\Pc$ and $q:Y\to\Pc$, let
\begin{equation}\label{E:tensorcoprod}
X\Times{p}{q}Y=\{(x,y)\in X\Tim Y : p(x)=q(y)\}\,.
\end{equation}
A groupoid is a small category in which all arrows are invertible. It consists of a set of
arrows $\G$, a set of objects $\Pc$ (called the {\em base}),
 source and end maps  ${\mathfrak s}, {\mathfrak e}: \G \to \Pc$,
 composition $m:\G\Times{\frake}{\fraks}\G\to\G$, and identities $\id:\Pc\to\G$.
We use $\G$ to denote both the groupoid and the set of arrows.
We write $m(f, g) = fg$ (and not $gf$) for $(f,g)\in\G\Times{\frake}{\fraks}\G$
 and we usually identify $\Pc$ with the subset $\id(\Pc)$ of $\G$.
For $P, Q\in \Pc$ we let ${}_P\G=\{x\in \G: \fraks(g) = P\}$, $\G_P= \{x\in \G: \frake(g) = P\}$,
$\G(P,Q)={}_P\G\cap \G_Q$, and $\G(P) = \G(P, P)$.

Alternatively, a groupoid may be defined as a set $\G$ with a partially defined associative
product and partial units, whose  elements are all invertible. The source and end maps
are determined by $\fraks(g) = gg^{-1}$, $\frake(g) = g^{-1}g$, and
 $P\in \Pc$ if and only if $P^2 = P$.

The opposite groupoid to $\G$ (where $\fraks$ and $\frake$ are switched) is denoted $\G^{op}$.

A  morphism of groupoids from $\G$
to $\K$ is a functor $\alpha: \G \to \K$.  Equivalently, $\alpha$ is a map from $\G$ to $\K$ which
preserves the product, and hence also
the base, source and end maps.  If $\G$ and $\K$ have the same base
$\Pc$, we say that $\alpha: \G \to \K$ is a morphism of groupoids \emph{over} $\Pc$ if
the restriction $\Pc \to \Pc$ is the identity.

The category of groupoids over  $\Pc$ has an initial object, the {\em discrete} groupoid
with $\fraks=\frake=m=\id$ are all the identity map $\Pc \to \Pc$, and a final one, the {\em coarse} groupoid $\Pc \times \Pc$,
where $\fraks(P,Q)=P$, $\frake(P,Q)=Q$, $\id(P)=(P,P)$ and 
$m\bigl((P,Q),(Q,R)\bigr)=(P,R)$.

A group bundle is a groupoid $\N$ with $\fraks=\frake$; thus $\N = \coprod_{P\in \Pc} \N(P)$.
Let $\alpha: \G \to \K$ be a morphism of groupoids over  $\Pc$. The kernel
of $\alpha$ is the group bundle
$\N = \{g\in \G: \alpha(g) \in \Pc\}$. In this case
 $\N(P) = \ker (\alpha: \G(P) \to \K(P))$.

\subsection{Actions and matched pairs of groupoids}

\

Let $\G$ be a groupoid with base $\Pc$.
Let $p: \Ec \to \Pc$ be a map. A \emph{left action} of $\G$ on $p$ is a map
$\fde : \G \Times{\mathfrak e}{p} \Ec \to \Ec$ such that
\begin{flalign}
\label{mp-0}  & p(g\fde e) = {\mathfrak s}(g),&
\\ \label{mp-1}  & g \fde(h \fde e)  = gh \fde e,&
\\ \label{mp-1.5}  & \id \,{p(e)} \, \fde  e  =  e,&
\end{flalign}
for all $g, h \in \G$, $e \in \Ec$  composable in the appropiate sense.
We may represent~\eqref{mp-0} by
\[\xymatrix{
{\Ec}\ar[d]_{p} & {g\fde e}\ar@{|->}[d]  & {e}\ar@{|->}[d]\\ 
{\Pc} & {Q}\ar@/^/[r]^{g} & {P}
}\]

Suppose $\G$ acts on
$p: \Ec \to \Pc$ and on $p': \Ec' \to \Pc$. A map $\psi:\Ec\to\Ec'$ {\em intertwines} the actions  if
\begin{equation}\label{E:intertwiner}
p'(\psi(e))=p(e) \text{ \ and \ } \psi(g\fde
e)=g\fde\psi(e)\,.
\end{equation}

There are two kinds of {\em trivial} actions. First, if $\Ec$ is of the form
$\Ec = \Pc\times X$ for some set $X$ and $p$ is the projection onto the first coordinate,
then an action on $p$ of an arbitrary groupoid $\G$ with base $\Pc$ is called trivial if 
 \begin{equation}\label{E:trivialuno}
g\fde ({\mathfrak e}(g), x) =({\mathfrak s}(g), x)
\end{equation}
 for all $x\in X$, $g\in \G$. Second, if $\N$ is a group bundle over $\Pc$,
 an action of $\N$ on an arbitrary map $p:\Ec\to\Pc$ is called trivial if
\begin{equation}\label{E:trivialdos}
n\fde e = e
\end{equation}
for all $(n,e)\in\N\Times{\frake}{p}\Ec$. If $\N$ is a group bundle and $\Ec = \Pc\times X$ then
the two kinds of trivial actions agree.

\medbreak
Similarly, a \emph{right action} of $\G$ on a map $q:\Ec\to\Pc$ is a map
$\fiz: \Ec \Times{q}{\mathfrak s} \G \to \Ec$ such that
\begin{flalign}
\label{mp-0.3}  &q(e\fiz g) = {\mathfrak e}(g),&
\\ \label{mp-2} & (e \fiz g) \fiz h  = e \fiz gh, &
\\ \label{mp-2.5} & e  \fiz \, \id \,{q(e)}  = e, &
\end{flalign}
for all $g, h \in \G$, $e \in \Ec$  composable in the appropiate sense.
Any left action gives rise to a right action on the same map by $e\fiz g = g^{-1}\fde e$, and vice versa.

\begin{definition}\label{matchpairgpds}  \cite[Definition 2.14]{mk1}.
A \emph{matched pair of groupoids} is a pair of
groupoids $(\Vc,\Hc)$ with the same base $\Pc$ together with the following data.
Let $t, b: \Vc\rightrightarrows\Pc$ be the source and end maps of $\Vc$, respectively,
and $l, r:\Hc\rightrightarrows\Pc$ the source and end maps of $\Hc$, respectively.
The data consists of a left action $\fde : \Hc \Times{r}{t} \Vc \to \Vc$ of $\Hc$ on
$t: \Vc \to \Pc$, and a right action $\fiz : \Hc \Times{r}{t} \Vc \to \Hc$ of $\Vc$ on
$r: \Hc \to \Pc$,
satisfying
\begin{flalign}
\label{mp-0.7}  &b(x\fde g) = l(x \fiz g),&
\\ \label{mp-3} & x \fde fg  = (x \fde f) ((x \fiz f) \fde g), &
\\ \label{mp-4} & xy \fiz g = (x \fiz (y \fde g)) (y \fiz g), &
\end{flalign}
for all $f, g \in \Vc$, $x, y \in \Hc$ for which the operations are defined.
\end{definition}

We refer to $\Vc$  as the {\em vertical} groupoid, $t$ and $b$ stand for {\em top} and {\em bottom};
$\Hc$ is the {\em horizontal} groupoid, $l$ and $r$ stand for {\em left} and {\em right}. The following
diagram illustrates the situation:
\begin{equation}\label{E:cell}
\xymatrix@C+5pt{
{t(x\fde g)=l(x)}\ar@/_/[d]_{x\fde g}\ar@/^1pc/[r]^{x} & {r(x)=t(g)}\ar@/^/[d]^{g}\\  
{b(x\fde g)=l(x\fiz g)}\ar@/_1pc/[r]_{x\fiz g} & {r(x\fiz g)=b(g)} }
\end{equation}
We refer to such diagrams as the {\em cells} of the matched pair.
We will not need any notions from double category theory beyond
this basic notation. For more in this direction, see~\cite{mk1,AN}.

Given $P\in \Pc$, there are two identities: $\id_{\Hc} P \in \Hc$
and $\id_{\Vc}P \in \Vc$. 

\begin{lema} For all $x,y \in \Hc$ and $f, g \in \Vc$ for which the operations are defined, we have
\begin{align}\label{iduno}
x \fde \id_\Vc r(x) &= \id_\Vc l(x)\,, \\
\label{iddos}
\id_\Hc t(g) \fiz g &= \id_\Hc b(g)\,,\\
\label{idtres}
(x \fde g)^{-1} &= (x \fiz g) \fde g^{-1}\,,\\
\label{idcuatro}
(x \fiz g)^{-1} &= x^{-1} \fiz (x \fde g)\,,\\
\label{idcinco}
(x \fiz g)^{-1} \fde (x \fde g)^{-1} &= g^{-1}\,,\\
\label{idseis}
(x \fiz g)^{-1} \fiz (x \fde g)^{-1} &= x^{-1}\,,\\
\label{idsiete}
g^{-1} (x^{-1} \fde f) &= (x \fiz g)^{-1} \fde \left((x \fde g)^{-1} f\right)\,, \\
\label{idocho}
(y \fiz g^{-1})x^{-1}  &= \left(y (x \fiz g)^{-1} \right) \fiz (x \fde g)^{-1}\,.
\end{align}
\end{lema}
\pf  First, by~\eqref{mp-1.5}, $x\fiz\id_{\Vc}r(x)=x$; hence, by~\eqref{mp-3},
\begin{multline*}
x\fde\id_{\Vc} r(x)=x\fde \bigl(\id_{\Vc}r(x)\id_{\Vc}r(x)\bigr) \\
=\bigl(x \fde \id_{\Vc}r(x)\bigr)
\Bigl(\bigl(x\fiz \id_{\Vc}r(x)\bigr) \fde \id_{\Vc}r(x)\Bigr)=\bigl(x \fde \id_{\Vc}r(x)\bigr)\bigl(x
\fde \id_{\Vc}r(x)\bigr)\,.
\end{multline*}
 Cancelling we obtain \eqref{iduno}. Now, 
\[\id_\Vc l(x) \equal{\eqref{iduno}} x\fde (gg^{-1}) \equal{\eqref{mp-3}} (x \fde g) \bigl((x \fiz
g) \fde g^{-1}\bigr)\,,\]
from which~\eqref{idtres} follows.  Furthermore,
$(x \fiz g)^{-1} \fde (x \fde g)^{-1} =(x \fiz g)^{-1} \fde \left( (x \fiz g) \fde g^{-1}\right) = 
g^{-1}$, which proves \eqref{idcinco}. Finally,
\[(x \fiz g)^{-1} \fde \left((x \fde g)^{-1} f\right)
\equal{\eqref{mp-3}} \left((x \fiz g)^{-1} \fde (x \fde g)^{-1}\right)\left(\left((x \fiz g)^{-1} \fiz
(x \fde g)^{-1}\right) \fde f\right) = g^{-1} (x^{-1} \fde f)\, \]
which proves~\eqref{idsiete}. The proofs of~\eqref{iddos},~\eqref{idcuatro},~\eqref{idseis},
and~\eqref{idocho} are similar.
\epf

Conditions~\eqref{idcinco} and~\eqref{idseis} may be succinctly expressed
as follows: inverting all arrows in the cell~\eqref{E:cell} yields
a new cell, as shown below:
\begin{equation}\label{E:inversecell}
\xymatrix@C+5pt{
{r(x\fiz g)=b(g)}\ar@/_/[d]_{g^{-1}}\ar@/^1pc/[r]^{(x\fiz g)^{-1}} & {b(x\fde g)=l(x\fiz
g)}\ar@/^/[d]^{(x\fde g)^{-1}}\\   {r(x)=t(g)}\ar@/_1pc/[r]_{x^{-1}} & {t(x\fde g)=l(x)} }
\end{equation}

\begin{obs}\label{R:basic} As for Hopf algebras (or weak Hopf algebras), there is a number
of basic constructions one may perform with matched pairs of groupoids. Let $(\Vc,\Hc)$
be a matched pair. The {\em dual} or {\em transpose} matched pair is $(\Hc,\Vc)$, with the
following actions:
\begin{equation}\label{E:dualactions}
f\rightharpoonup y:= (y^{-1}\fiz f^{-1})^{-1} \text{ \ and \ }f\leftharpoonup y:= (y^{-1}\fde
f^{-1})^{-1}
\end{equation}
for $(f,y)\in\Vc\Times{b}{l}\Hc$. Note that if $(x,g)\in\Hc\Times{r}{t}\Vc$ then
$(x\fde g,x\fiz g)\in\Vc\Times{b}{l}\Hc$. According  to~\eqref{idcinco} and~\eqref{idseis}, the
corresponding cell of the matched pair $(\Hc,\Vc)$ is
\begin{equation}\label{E:dualcell}
\xymatrix@C+5pt{
{P}\ar@/_/[d]_{x}\ar@/^/[r]^{x\fde g} & {Q}\ar@/^/[d]^{x\fiz g}\\  
{R}\ar@/_/[r]_{g} & {S} }
\end{equation}
In other words, transposing a cell~\eqref{E:cell} of the matched
pair $(\Vc,\Hc)$ yields a cell of the dual matched pair $(\Hc,\Vc)$.

The dual of the dual of $(\Vc,\Hc)$ is the original matched pair $(\Vc,\Hc)$.

The {\em opposite} of $(\Vc,\Hc)$  is the matched pair $(\Vc,\Hc^{op})$ with the following actions:
\begin{equation}\label{E:oppactions}
y\fde^{op} g:= y^{-1}\fde g \text{ \ and \ } y\fiz^{op} g:= (y\fiz g)^{-1}
\end{equation}
for $(f,y)\in(\Hc^{op})\Times{{r^{op}}}{t}\Vc=\Hc\Times{l}{t}\Vc$. Note that if
$(x,g)\in\Hc\Times{r}{t}\Vc$ then $(x^{-1},g)\in\Hc\Times{l}{t}\Vc$. The
corresponding cell of the matched pair $(\Vc,\Hc^{op})$ is
\begin{equation}\label{E:oppcell}
\xymatrix@C+5pt{
{Q}\ar@/_/[d]_{x\fde g}\ar@/^/[r]^{x^{-1}} & {P}\ar@/^/[d]^{g}\\  
{S}\ar@/_/[r]_{(x\fiz g)^{-1}} & {R} }
\end{equation}
In other words, inverting the horizontal arrows in a cell~\eqref{E:cell} of the matched
pair $(\Vc,\Hc)$ yields a cell of the opposite matched pair $(\Hc,\Vc)$.

The {\em coopposite} matched pair $(\Vc^{op},\Hc)$ is defined similarly, by inverting the
vertical arrows.
\end{obs}

\begin{exa}\label{Ex:initial} Let $\Hc=\Pc$ be the discrete groupoid over $\Pc$ and $\Vc=\Pc\times\Pc$ be
the coarse groupoid over $\Pc$ (Section~\ref{S:basic}). Then $(\Vc,\Hc)$ is a matched pair of groupoids
with the following actions:
\[P\fde (P,Q)=(P,Q) \text{ \ and \ }P\fiz(P,Q)=Q\,.\]
The cells of this matched pair~\eqref{E:cell} are in this case simply
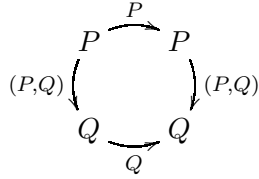
\[\xymatrix@-5pt{
{P}\ar@/_/[d]_{(P,Q)}\ar@/^/[r]^{P} & {P}\ar@/^/[d]^{(P,Q)}\\  
{Q}\ar@/_/[r]_{Q} & {Q} }\]

The dual matched pair has $\Hc=\Pc\times\Pc$, $\Vc=\Pc$, and actions
\[(P,Q)\fde Q = P \text{ \ and \ } (P,Q)\fiz Q = (P,Q)\,. \]
We refer to $(\Pc\times\Pc,\Pc)$ as the {\em initial} matched pair and to $(\Pc,\Pc\times\Pc)$
as the {\em terminal} matched pair. This terminology is justified by Proposition~\ref{P:initial}.
\end{exa}

\begin{exa}\label{Ex:XY} Let $X$ and $Y$ be arbitrary sets. Define $\Pc:=X\times Y$,
$\Hc:=X\times X\times Y$, and $\Vc:=X\times Y\times Y$. In order to define groupoid
structures on $\Hc$ and $\Vc$ with base $\Pc$ and a matched pair structure on $(\Vc,\Hc)$, 
we first specify the cells:
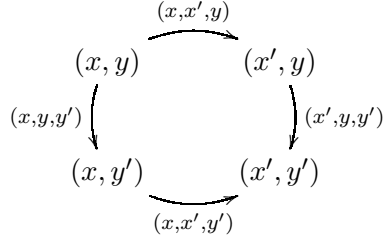
\[\xymatrix@C+5pt{
{(x,y)}\ar@/_/[d]_{(x,y,y')}\ar@/^1pc/[r]^{(x,x',y)} & {(x',y)}\ar@/^/[d]^{(x',y,y')}\\  
{(x,y')}\ar@/_1pc/[r]_{(x,x',y')} & {(x',y')} }\]
{}From this diagram we can tell the source and end maps of $\Hc$ and $\Vc$, as well as the
actions $\fde$ and $\fiz$. For instance,
\[l(x,x',y):=(x,y)\,,\  b(x,y,y'):=(x,y')\,, \text{ and } (x,x',y)\fde (x',y,y'):=(x,y,y')\,.\]
It remains to describe the products on $\Hc$ and $\Vc$. They are
\[m_{\Hc}\bigl((x,x',y),(x',x'',y)\bigr)=(x,x'',y) \text{ \ and \ }
m_{\Vc}\bigl((x,y,y'),(x,y',y'')\bigr)=(x,y,y'')\,. \]
The groupoid and matched pair axioms are easily verified. We denote this matched pair by
$M(X,Y)$. Note that if $X$ is a singleton then $M(X,Y)=(\Pc\times\Pc,\Pc)$, the initial matched
pair, while if $Y$ is a singleton then $M(X,Y)=(\Pc,\Pc\times\Pc)$, the terminal matched
pair. The dual of $M(X,Y)$ is $M(Y,X)$.

For a characterization of matched pairs of the form $M(X,Y)$ see~\cite[Proposition 2.14]{AN}.
\end{exa}

\begin{exa}\label{Ex:semi-initial} Let $\Vc$ be any groupoid with base $\Pc$. There is a matched
pair $(\Vc,\Pc)$ with actions
\[t(f)\fde f=f \text{ \ and \ }t(f)\fiz f=b(f)\,.\]
Similarly, for any groupoid $\Hc$ with base $\Pc$, there is a matched pair
$(\Pc,\Hc)$  with actions
\[x\fde r(x) = l(x) \text{ \ and \ } x\fiz r(x) = x\,. \]
The dual of $(\Vc,\Pc)$ is $(\Pc,\Vc)$ (and the dual of $(\Pc,\Hc)$ is $(\Hc,\Pc)$).
\end{exa}

\begin{exa}\label{Ex:automorphisms}  Let $\Vc$ be a groupoid with base $\Pc$, source $t$ and end $b$. 
Let $\N$ be a group bundle over $\Pc$, viewed as a groupoid with 
source and end $l=r$. Consider the trivial left
action of
$\N$ on $t:\Vc\to\Pc$ as in~\eqref{E:trivialdos}: $n\fde g = g$ for all $(n,g) \in \N\Times{r}{t}\Vc$.

Consider an arbitrary right action of $\Vc$ on $r:\N\to\Pc$. Then axioms~\eqref{mp-0.7} and~\eqref{mp-3}
are satisfied, while axiom~\eqref{mp-4} is satisfied if and only if
\[nm\fiz g = (n\fiz g)(n\fiz g)\] 
for all composable $n,m\in \N$, $g \in\Vc$.
If this is the case, we say that the action of $\Vc$ on $\N$ is \emph{by group
bundle automorphisms}. Thus, an action of a groupoid $\Vc$ on a group bundle $\N$ by automorphisms 
together with the
trivial action of the group bundle on the groupoid yield a matched pair of groupoids
$(\Vc,\N)$. 

There is a similar notion of left action by group bundle automorphisms
of a a groupoid $\Hc$ on a group bundle $\N$. Together with the trivial right
action of $\N$ on $\Hc$, this gives rise to a matched pair $(\N,\Hc)$.

There is a special case of particular interest. Let $\Hc$ be an arbitrary
groupoid with base $\Pc$ and let $\N:=\coprod_{P\in\Pc}\Hc(P,P)$. 
Then $\Hc$ acts on $\N$ by conjugation: $x\fde n:=xnx^{-1}$. The matched pair
$(\N,\Hc)$ is one of the simplest instances of the double construction 
(see Example~\ref{Ex:easydouble}).
\end{exa}

\subsection{Two groupoids associated to a matched pair}

\

\begin{definition}\label{D:diagonal} Let $(\Vc,\Hc)$ be a matched pair of groupoids with base $\Pc$. The
{\em diagonal} groupoid $\Vc\bowtie \Hc$ is defined as follows. The base is $\Pc$, the set of arrows is
$\Vc\bowtie \Hc := \Vc \Times{b}{l} \Hc$, the source map is
${\mathfrak s}(f, y) = t(f)$, the end map is
${\mathfrak e}(f, y) = r(y)$, the identities are $\id P= (\id_{\Vc} P,\id_{\Hc} P)$,
 and the  composition is
\begin{equation}\label{E:diag-product}
(f,y) (h, z) = (f (y\fde h), (y\fiz h) z)\,.
\end{equation}
\end{definition}

The elements $(f,y)$ of $\Vc\bowtie\Hc$ may be represented by {\em corner diagrams}
\[\xymatrix{
{t(f)}\ar@/_/[d]_{f} & \\  
{b(f)=l(y)}\ar@/_1pc/[r]_{y} & {r(y)} }\]
In this notation, the composition of $\Vc\bowtie\Hc$ results from the matched pair structure by stacking
two such corners diagonally and filling in with a cell:
\[\xymatrix{
{t(f)}\ar@/_/[d]_{f} & & \\  
{b(f)=l(y)}\ar@/^1pc/[r]^{y}\ar@{-->}@/_/[d]_{y\fde h} & {r(y)=t(h)}\ar@/^/[d]^{h} & \\
{b(y\fde h)=l(y\fiz h)}\ar@{-->}@/_1pc/[r]_{y\fiz h} & {b(h)=l(z)}\ar@/_1pc/[r]_{z} & {r(z)}
}\]

Note that, in view of~\eqref{iduno} and~\eqref{iddos}, $\Vc$ and $\Hc$ may be viewed as subgroupoids of
$\Vc\bowtie\Hc$ via
$f\mapsto(f,\id_\Hc b(f))$ and $y\mapsto(\id_\Vc l(y),y)$.

Let $\D$ be a groupoid with base $\Pc$. An \emph{exact factorization} of
$\D$ is a pair of subgroupoids $\Vc$, $\Hc$ with the  same base $\Pc$, such that for any $\alpha \in \D$,
there exist unique $f\in \Vc$ and $y \in \Hc$ with $\alpha = fy$; in other words, the composition
map $\Vc \Times{\frake}{\fraks} \Hc \to \D$ is a bijection. A matched pair of groupoids $(\Vc,\Hc)$ gives
rise to an exact factorization of the diagonal groupoid $\D:=\Vc\bowtie\Hc$ with factors $\Vc$ and $\Hc$.

Conversely, an exact factorization gives rise to a matched pair. Indeed, let $t$ and $b$ denote the
restrictions of $\fraks$ and $\frake$ to $\Vc$, and $l$ and $r$ the restrictions of $\fraks$ and $\frake$
to $\Hc$. Then $\Hc \Times{r}{t} \Vc\subseteq \D\Times{\frake}{\fraks}\D$. Now, given  $(x,g) \in \Hc
\Times{r}{t} \Vc$, the bijectivity of
the composition allows to define elements $x\fde g\in \Vc$ and $x\fiz g\in\Hc$ by the equation
\[xg = (x\fde g)(x\fiz g)\,.\]
Then $\Vc$ and $\Hc$, together with these actions, form a matched pair. 

The above constructions of exact factorizations and matched pairs are inverse. 

For more details on these constructions, the reader may consult~\cite[Theorems 2.10 and 2.15]{mk1}
or~\cite[Prop 2.9]{AN}.

\begin{exa} Consider a right action of a groupoid $\Vc$ on a group bundle $\N$
by automorphisms, and the corresponding matched pair $(\Vc,\N)$ 
(Example~\ref{Ex:automorphisms}).
In this case, we  denote
the diagonal groupoid by $\Vc \ltimes \N$ and refer to it as the semidirect 
product of $\Vc$ by $\N$.

Note that the projection $\Vc \ltimes \N \to \Vc$, $(g,n)\mapsto g$ 
is a morphism of groupoids,
and it has a section $\sigma: \Vc \to \Vc \ltimes \N$, 
$\sigma(g) = (g, \id b(g))$, which is a morphism of
groupoids. 

Conversely let $\alpha: \G \to \K$ be a morphism of groupoids over $\Pc$.  
If there is a  section $\sigma: \K \to \G$ that is a 
morphism of groupoids, then $\K$ acts on the kernel $\N$ of $\alpha$ by 
$n\fiz k = \sigma(k)^{-1}n\sigma(k)$ and
$\G \simeq \K \ltimes \N$. This may be seen as a special case of the
characterization of matched pairs in terms of exact factorizations.
\end{exa}

A simpler construction of groupoids is the following.

\begin{definition}\label{dirprodgpds}
Let $\G$ and $\K$ be groupoids over $\Pc$. The \emph{restricted product} of $\G$ and $\K$ is
$$
\G \prode \K := (\G\Times{\fraks}{\fraks}\K)\cap(\G\Times{\frake}{\frake}\K)=
\{(g, k) \in  \G \times \K: {\mathfrak s}(g) = {\mathfrak
s}(k),\ {\mathfrak e}(g) = {\mathfrak e} (k)\}.
$$
Then $\G \prode \K$ is a groupoid over the same base $\Pc$ with source and end maps defined by
\[\fraks(g,k)=\fraks(g)=\fraks(k) \text{ \ and \ } \frake(g,k)=\frake(g)=\frake(k) \,,\]
and composition
\[(g,k)(g',k')=(gg',kk')\,.\]
\end{definition}

The restricted product is the product in the category of groupoids over $\Pc$: given a groupoid $\F$
with base $\Pc$ and morphisms $\alpha:\F\to\G$ and $\beta:\F\to\K$ of groupoids over $\Pc$, there is a
unique morphism 
$\F\to \G \prode \K$ of groupoids over $\Pc$ fitting in the commutative diagram
\begin{equation}\label{E:universal-prode}
\xymatrix{
& {\F}\ar[ld]_{\alpha}\ar@{-->}[d]\ar[rd]^{\beta}\\
{\G} & {\G \prode \K}\ar[l]\ar[r] & {\K}
}
\end{equation}

\medbreak
Let $(\Vc,\Hc)$ be a matched pair.
The diagonal groupoid, on the other hand, may be viewed as a sort of ``restricted coproduct'',
in the sense that it satisfies the following universal property. Given a groupoid $\F$ with base $\Pc$ 
and morphisms $\alpha:\Hc\to\F$ and $\beta:\Vc\to\F$ of groupoids over $\Pc$ such that
\begin{equation}\label{E:commuting-diagonal}
\alpha(y)\beta(h)=\beta(y\fde h)\alpha(y\fiz h) \text{ \ for
every \ }(y,h)\in\Hc\Times{r}{t}\Vc\,,
\end{equation}
there is a unique morphism $\Vc\bowtie\Hc\to\F$ of groupoids over
$\Pc$ fitting in the commutative diagram
\begin{equation}\label{E:universal-diagonal}
\xymatrix{
{\Hc}\ar[r]\ar[rd]_{\alpha} & {\Vc\bowtie\Hc}\ar@{-->}[d] & {\Vc}\ar[l]\ar[ld]^{\beta}\\
&{\F}&\\
}
\end{equation}

Condition~\eqref{E:commuting-diagonal} may be illustrated as follows: 
\[\xymatrix@C+10pt{
{P}\ar@/_/[d]_{y\fde h}\ar@/^/[r]^{y} & {Q}\ar@/^/[d]^{h}\\  
{R}\ar@/_/[r]_{y\fiz h} & {S} }
\qquad \raisebox{-19pt}{$\Rightarrow$} \qquad
\raisebox{13pt}{\xymatrix@R-10pt{
& {Q}\ar[rd]^{\beta(h)} &\\
{P}\ar[ru]^{\alpha(y)}\ar[rd]_{\beta(y\fde h)} & & {S}\\
&{R}\ar[ru]_{\alpha(y\fiz h)} &} }
\]
In words, cells in the matched pair $(\Vc,\Hc)$
are transformed by $(\alpha,\beta)$ into {\em commutative} diagrams in the groupoid $\F$. 

\subsection{Morphisms of matched pairs. Duals and doubles}\label{S:morphisms}

\ 

\begin{definition}\label{D:morphism}
Let  $(\Vb,\Hb)$ and $(\Vc,\Hc)$ be two matched pairs of 
groupoids with the same base $\Pc$. A \emph{morphism of matched pairs of groupoids} 
$(\Vb,\Hb)\to (\Vc,\Hc)$  is a pair $(\alpha,\beta)$ of morphisms
of groupoids over $\Pc$, $\alpha: \Hb \to \Hc$ and $\beta: \Vc \to \Vb$,  such that
\begin{align}
\label{functor1}
\beta (\alpha (h) \fde g) &= h \fde \beta(g)\,,  \\
\label{functor2}
\alpha (h \fiz \beta (g)) &= \alpha (h) \fiz g\,,
\end{align}
for all $(h,g)\in \Hb\Times{r}{t}\Vc$.
\end{definition}

If $(\Hg,\Vg)$ is another matched pair of groupoids with  base $\Pc$
and $(\delta, \omega):(\Vc,\Hc)\to (\Vg,\Hg)$ is a morphism of matched pairs,
then $(\alpha\delta, \omega\beta)$ is a morphism of matched pairs  $(\Vb,\Hb)\to (\Vg,\Hg)$.
Thus, matched pairs of groupoids with the same base $\Pc$ form a category.

There are two distinguished objects in this category.

\begin{proposition}\label{P:initial} The matched pairs $(\Pc\times\Pc,\Pc)$ and  $(\Pc,\Pc\times\Pc)$ of
Example~\ref{Ex:initial} are respectively 
the initial and the terminal objects in the category of matched pair of groupoids with base $\Pc$. 
\end{proposition}
\pf Given an arbitrary matched pair $(\Vc,\Hc)$, consider the unique morphisms of groupoids over $\Pc$,
$\alpha_0:\Pc\to\Hc$ and $\beta_0:\Vc\to\Pc\times\Pc$ (Section~\ref{S:basic}):
\[\alpha_0(P)=\id_{\Hc}P\,,\ \ \ \beta_0(g)=(t(g),b(g))\,. \]
Then $(\alpha_0,\beta_0)$ is (the unique) morphism of
matched pairs $(\Pc\times\Pc,\Pc)\to (\Vc,\Hc)$: condition~\eqref{functor1} holds by~\eqref{mp-1.5}
and condition~\eqref{functor2} holds by~\eqref{iddos}. Similarly, the unique morphisms of
groupoids over $\Pc$, $\alpha_1:\Hc\to\Pc\times\Pc$ and $\beta_1:\Pc\to\Vc$ combine to give the
unique morphism of matched pairs $(\alpha_1,\beta_1):(\Vc,\Hc)\to(\Pc,\Pc\times\Pc)$.
\epf

\begin{obs}\label{R:iso-opp} Any matched pair $(\Vc,\Hc)$ is isomorphic to its opposite $(\Vc,\Hc^{op})$
(Remark~\ref{R:basic}). The isomorphism consist of the pair
$\alpha:\Hc\to\Hc^{op}$, $\alpha(x)=x^{-1}$, and $\beta:\Vc\to\Vc$, $\beta(g)=g$. 
Similarly,  any matched pair is isomorphic to its coopposite.
\end{obs}

\bigbreak
We will construct a new matched pair from a morphism of matched pairs. The following lemma
is a preliminary step towards this goal. We use $t$ and $b$ to denote the source and end maps
of both $\Vb$ and $\Vc$, and similarly $l$ and $r$ for the source and end maps of both $\Hb$ and $\Hc$.

\begin{lema}\label{double-lema}  Let $(\alpha, \beta):(\Vb,\Hb)\to (\Vc,\Hc)$ be a morphism of
matched pairs of groupoids. Define maps $\rightharpoonup:\Hb \Times{r}{t} \Vc\to\Vc$ and 
$\leftharpoonup:\Hb \Times{r}{t} \Vc\to\Hb$ by
\begin{equation}\label{action10}
h \rightharpoonup g := \alpha(h) \fde g\,, \qquad
h \leftharpoonup g := h \fiz \beta(g)\,. 
\end{equation}
 Then $(\Vc,\Hb,\rightharpoonup,\leftharpoonup)$ is a matched pair
of groupoids.
\end{lema}

\pf The maps $\rightharpoonup$ and  $\leftharpoonup$ are well-defined and are actions
since $\alpha$ and $\beta$ are morphisms of groupoids. We check \eqref{mp-0.7}:
$b(h \rightharpoonup g) = b(\alpha(h) \fde g) = l(\alpha(h) \fiz g) =
l(\alpha(h \fiz \beta(g)) ) =  l (h \fiz \beta(g)) = l(h \leftharpoonup g)$,
where we have used \eqref{functor2} in the third equality.
We check \eqref{mp-3}:
$h \rightharpoonup fg = \alpha(h) \fde fg = (\alpha(h) \fde f)((\alpha(h) \fiz f) \fde g) =
(h \rightharpoonup f)(\alpha(h \leftharpoonup f) \fde g) =
(h \rightharpoonup f)( (h \leftharpoonup f) \rightharpoonup g) $,
where we have used \eqref{functor2} in the third equality.
We check \eqref{mp-4}:
$hk \leftharpoonup g = hk \fiz \beta(g) = (h \fiz ( k \fde \beta(g))) ( k \fiz \beta(g))=
(h \fiz \beta( k \rightharpoonup g)) ( k \leftharpoonup g) =
(h \leftharpoonup ( k \rightharpoonup g)) ( k \leftharpoonup g)$,
where we have used \eqref{functor1} in the third equality. \epf

We may now consider the diagonal groupoid $ \Vc \bowtie \Hb$ (Definition~\ref{D:diagonal})
associated to the matched pair of Lemma~\ref{double-lema}. Let $l$ and $r$ be its source and
end maps. According to Definition~\ref{D:diagonal},
\[l(g,h) = t(g) \text{ \ and \ }r(g,h) = r(h)\]
 for all $(g,h)\in \Vc\bowtie \Hb$; and the product is given by
\begin{equation}\label{E:diag-product-2}
(f,u)(g,h) = \Bigl(f\bigl(\alpha(u)\fde g\bigr), \bigl(u\fiz \beta(g)\bigr)h\Bigr)\,.
\end{equation}
The groupoid $\Vc\bowtie \Hb$ will be the horizontal groupoid of a matched pair to be defined shortly.
The vertical groupoid will be $\Vb \prode \Hc^{op}$. Let $t$ and $b$ be the source and
end maps of $\Vb \prode \Hc^{op}$. According to Definition~\ref{dirprodgpds}, 
 \[t(\gamma, x) =t(\gamma) = r(x) \text{ \ and \ }
b(\gamma, x) = b(\gamma) = l(x)\]
for all $(\gamma, x) \in\Vb\prode\Hc^{op}$. 

\begin{theorem}\label{T:double} 
 Let $(\alpha, \beta):(\Vb,\Hb)\to (\Vc,\Hc)$ be a morphism of
matched pairs of groupoids. Define maps $\rightharpoonup:(\Vc \bowtie \Hb) \Times{r}{t} (\Vb \prode
\Hc^{op})\to \Vb \prode \Hc^{op}$ and
$\leftharpoonup:(\Vc \bowtie \Hb) \Times{r}{t} (\Vb \prode \Hc^{op})
\to \Vc \bowtie \Hb$ by
\begin{align}\label{action11}
(g,h) \rightharpoonup (\gamma, x) &\ :=\ 
\Biggl(\beta(g) (h\fde \gamma)  \beta\Bigl(\bigl(\alpha(h\fiz \gamma)x\alpha(h)^{-1}\bigr) \fde
g^{-1}\Bigr)\,,\ 
\Bigl(\alpha(h\fiz \gamma)x\alpha(h)^{-1}\Bigr) \fiz g^{-1} \Biggr)\,,
\\\label{action11bis}
(g,h) \leftharpoonup (\gamma, x) &\ :=\   
\Biggl( \Bigl(\bigl(\alpha(h\fiz \gamma)x\alpha(h)^{-1}\bigr) \fiz g^{-1}\Bigr) \fde g\,,\  h\fiz
\gamma\Biggr)\,,
\end{align}
for all $(g,h)\in \Vc \bowtie \Hb$, $(\gamma, x) \in \Vb \prode \Hc^{op}$ with $r(h) = t(\gamma)$.

\bigbreak
Then $(\Vb \prode \Hc^{op},\Vc \bowtie \Hb, \rightharpoonup,\leftharpoonup)$ is a matched pair of
groupoids.
\end{theorem}

\medbreak

The proof of this theorem is lengthy but for the most part straightforward
and is given in the Appendix. Here we content ourselves with an
illustration of formulas~\eqref{action11} and~\eqref{action11bis}, in the
special case when $(\Vb,\Hb)=(\Vc,\Hc)$ and $\alpha$ and $\beta$ are identities.
We start from $(g,h)\in\Vc\bowtie\Hc$ and $(\gamma,x)\in\Vc\prode\Hc^{op}$
with $r(h)=t(\gamma)$ (so $b(g)=l(h)$, $t(\gamma)=r(x)$, and $b(\gamma)=l(x)$),
and fill in with cells:
\begin{equation}\label{E:prod-double}
\xymatrix@R-10pt@C+5pt{
 {P}\ar@/_/[d]_{g}\\
 {Q}\ar@/_/[d]_{h\fde\gamma}\ar@/^/[r]^{h} & {R}\ar@/^/[d]^{\gamma}\\
{R'}\ar@/_/[r]_{h\fiz \gamma}\ar@/_/[d]_{\Omega\fde g^{-1}} &{S}\ar@/_/[r]_{x} & {R}\ar@/_/[r]_{h^{-1}} & 
{Q}\ar@/^/[d]^{g^{-1}}\\
{Q'}\ar@/_1pc/[rrr]_{\Omega\fiz g^{-1}} & & & {P}  }
\end{equation}
where $\Omega:=(h\fiz\gamma)xh^{-1}\in\Hc$. Formulas~\eqref{action11} 
and~\eqref{action11bis} can be written as
\[(g,h) \rightharpoonup (\gamma, x) =
\bigl(g(h\fde \gamma)(\Omega \fde g^{-1})\,, 
\Omega \fiz g^{-1} \bigr) \text{ \ and \ }
(g,h) \leftharpoonup (\gamma, x) =   
\bigl( (\Omega \fde g^{-1})^{-1}\,,  h\fiz\gamma\bigr)\,.\]

\bigbreak

\begin{definition}\label{D:double}
The matched pair associated to the identity morphism $(\Vc,\Hc)\to(\Vc,\Hc)$ by
means of Theorem~\ref{T:double} is called the {\em double} of $(\Vc,\Hc)$ and denoted
$D(\Vc,\Hc)$. The matched pair
associated to a general morphism
$(\alpha, \beta):(\Vb,\Hb)\to (\Vc,\Hc)$ is called a {\em generalized double} and denoted
$D(\alpha, \beta)$. 
\end{definition}

This terminology is justified by Theorem~\ref{T:centralizer}, 
Corollary~\ref{C:center}, and Theorem~\ref{T:drinfeld} below. 

\begin{obs}\label{R:dual} Consider the unique morphism
$(\alpha_0,\beta_0):(\Pc\times\Pc,\Pc)\to (\Vc,\Hc)$ (Proposition~\ref{P:initial}).
We describe its generalized double. We have
\[(\Pc\times\Pc)\prode\Hc^{op}\cong\Hc^{op} \text{ \ and \ } \Vc\bowtie\Pc\cong \Vc\]
 via $\bigl( r(x),l(x),x\bigr)\leftrightarrow x$ and $\bigl(g, b(g)\bigr)\leftrightarrow g$.
After these identifications,~\eqref{action11} and~\eqref{action11bis} boil down to
\[g\rightharpoonup x = x\fiz g^{-1} \text{ \ and \ }g\leftharpoonup x = (x\fiz g^{-1})\fde
g=(x\fde g^{-1})^{-1}\]
for $g\in\Vc$ and $x\in\Hc$ with $b(g)=r(x)$.
Comparing these expressions with~\eqref{E:dualactions} and~\eqref{E:oppactions} we see that
the generalized double $D(\alpha_0,\beta_0)$ is the dual of the opposite of $(\Vc,\Hc)$.
Since any matched pair is isomorphic to its opposite (Remark~\ref{R:iso-opp}),
this generalized double is isomorphic to the dual of
$(\Vc,\Hc)$: $D(\alpha_0,\beta_0)\cong(\Hc,\Vc)$.

Consider the generalized double of the (unique) morphism $(\alpha_1,\beta_1):(\Vb,\Hb)\to
(\Pc,\Pc\times\Pc)$ of Proposition~\ref{P:initial}. We have
\[\Vb\prode(\Pc\times\Pc)^{op}\cong\Vb \text{ \ and \ } \Pc\bowtie\Hb\cong \Hb\,,\]
and after these identifications,~\eqref{action11} and~\eqref{action11bis} boil down to
the original actions of the matched pair $(\Vb,\Hb)$. Thus
$D(\alpha_1,\beta_1)\cong(\Vb,\Hb)$. 
\end{obs}

\begin{exa} Consider the initial matched pair
$(\Pc\times\Pc,\Pc)$  (Example~\ref{Ex:initial}). The only endomorphism of $(\Pc\times\Pc,\Pc)$ is the
identity, so the dual and the double of this matched pair coincide. Therefore, the double of this matched
pair is $(\Pc,\Pc\times\Pc)$, the terminal matched pair (Example~\ref{Ex:initial}).

The unique morphism $(\alpha_1,\beta_1):(\Pc,\Pc\times\Pc)\to
(\Pc,\Pc\times\Pc)$ is the identity. Therefore,
by Remark~\ref{R:dual}, $(\Pc, \Pc\times\Pc)$ coincides with its own double.
\end{exa}

\begin{exa}\label{Ex:XYdouble} Let $X$ and $Y$ be sets and consider the matched pair $M(X,Y)$
of Example~\ref{Ex:XY}. We have $\Pc=X\times Y$, $\Hc=X\times X\times Y$, $\Vc=X\times Y\times Y$, and
\[\Vc\prode\Hc^{op}\cong X\times Y \text{ \ and \ } \Vc\bowtie\Hc\cong (X\times Y)\times(X\times Y)\]
via $\bigl((x,y,y),(x,x,y)\bigr)\leftrightarrow (x,y)$ and
$\bigl((x,y,y''),(x,x'',y'')\bigr)\leftrightarrow \bigl((x,y),(x'',y'')\bigr)$. Thus, the double
of this matched pair is the terminal matched pair $M(X\times Y,\{\ast\})=(\Pc,\Pc\times\Pc)$.
\end{exa}

\begin{exa}\label{Ex:easydouble}
Let $(\Vc,\Pc)$ and $(\Pc,\Hc)$ be the matched pairs of 
Example~\ref{Ex:semi-initial}. For $(\Pc,\Hc)$, we have
 $\Pc\prode\Hc^{op}=\coprod_{P\in\Pc}\Hc(P,P)^{op}$ (a group bundle)
and $\Pc\bowtie\Hc\cong\Hc$. Actions~\eqref{action11} and~\eqref{action11bis}
reduce to
\[h\rightharpoonup x= hx h^{-1} \text{ \ and \ }
h\leftharpoonup x =h\]
for $h\in\Hc(P,Q)$, $x\in\Hc(Q,Q)$. Thus, the double of $(\Pc,\Hc)$ is the coopposite of
(and hence isomorphic to) the matched pair 
$(\coprod_{P\in\Pc}\Hc(P,P),\Hc)$ of Example~\ref{Ex:automorphisms}. 

Similarly, the double of $(\Vc,\Pc)$ is the matched pair
$(\coprod_{P\in\Pc}\Vc(P,P),\Vc)$ of Example~\ref{Ex:automorphisms}.
In particular, $D(\Vc,\Pc)\cong D(\Pc,\Vc)$.
\end{exa}

Next we discuss the functoriality of generalized doubles.

\begin{proposition}\label{P:functoriality} 
Consider a commutative diagram of morphisms of matched pairs: 
\[\xymatrix@R-1pc@C-1.5pc{
{(\Vb,\Hb)}\ar[rd]_{(\varphi,\psi)}\ar[rr]^{(\alpha,\beta)} & &
{(\Vc,\Hc)}\ar[ld]^{(\chi,\omega)}\\ &{(\Vg,\Hg)} }\] 
There are morphisms of matched pairs
\[(\alpha_\#,\beta_\#):D(\varphi,\psi)\to D(\chi,\omega)
\text{ \ and \ } (\chi^\#,\omega^\#):D(\varphi,\psi)\to D(\alpha,\beta)\]
given by
\[\alpha_\#:=\Vg\bowtie\alpha\,,\ \beta_\#:=\beta\prode\Hg^{op}\,,\ 
\chi^\#:=\chi\bowtie\Hb\,,\ \omega^\#:=\Vb\prode\omega\,.\]
\end{proposition}
\pf We only note that $D(\varphi,\psi)=(\Vb\prode\Hg^{op},\Vg\bowtie\Hb)$,
$D(\chi,\omega)=(\Vc\prode\Hg^{op},\Vg\bowtie\Hc)$ and
$D(\alpha,\beta)=(\Vb\prode\Hc^{op},\Vc\bowtie\Hb)$, so the above maps
are well-defined, and omit the remaining verifications.
\epf

\begin{obs}\label{R:inclusion}  
Let us apply Proposition~\ref{P:functoriality} to the
commutative diagram
\[\xymatrix@R-1pc@C-1.5pc{
{(\Vb,\Hb)}\ar[rd]_{}\ar[rr]^{(\alpha,\beta)} & &
{(\Vc,\Hc)}\ar[ld]^{}\\ &{(\Pc,\Pc\times\Pc)} }\] 
where the diagonal maps are the unique morphisms to the
terminal matched pair. According to Remark~\ref{R:dual},
the generalized double of the left diagonal map can be identified with $(\Vb,\Hb)$. 
Therefore, the right diagonal map induces a morphism of matched pairs
\[(\Vb,\Hb)\to D(\alpha,\beta)\,.\]
Explicitly, this consists of
 the pair of morphisms of groupoids $\iota:\Hb\to \Vc \bowtie \Hb$
and $\pi:\Vb \prode \Hc^{op}\to\Vb$ given by
\[\iota(h)=(\id_{\Vc}l(h),h) \text{ \ and \ }\pi(\gamma,x)=\gamma\,.\]
This morphism plays a crucial role in Section~\ref{S:braidings}.
\end{obs}

Finally, let us apply Proposition~\ref{P:functoriality} to the
commutative diagram
\[\xymatrix@R-1pc@C-1.5pc{
{(\Vb,\Hb)}\ar[rr]^{(\alpha,\beta)} & &{(\Vc,\Hc)}\\ 
&{(\Pc\times\Pc,\Pc)}\ar[ru]\ar[lu] }\] 
where the diagonal maps are the unique morphisms from the
initial matched pair. By Remark~\ref{R:dual}, we obtain a morphism of matched pairs
$(\alpha^\#,\beta^\#)$ from the dual of $(\Vc,\Hc)$ to the dual of $(\Vb,\Hb)$.
Explicitly, $\alpha^\#=\beta:\Vc\to\Vb$ and $\beta^\#=\alpha:\Hb\to\Hc$.
We refer to the morphism $(\beta,\alpha)$ as the dual of the 
morphism $(\alpha,\beta)$.

The generalized double of a morphism and that of its dual coincide:
\begin{proposition}\label{P:gendouble-dual}
Let $(\alpha, \beta):(\Vb,\Hb)\to (\Vc,\Hc)$ be a morphism of matched
pairs and consider its dual $(\beta,\alpha):(\Hc,\Vc)\to(\Hb,\Vb)$. There is
an isomorphism of matched pairs
\[D(\alpha,\beta)\cong D(\beta,\alpha)\,.\]
\end{proposition}
\pf We define the isomorphism and omit the verifications. We have
\[D(\alpha,\beta)=(\Vb \prode \Hc^{op},\Vc \bowtie \Hb)
\text{ \ and \ }
D(\beta,\alpha)=(\Hc\prode\Vb^{op},\Hb\bowtie\Vc)\,.\]
The isomorphism $D(\alpha,\beta)\to D(\beta,\alpha)$ is given by the following pair of morphisms of groupoids:
\begin{gather*}
\Vc \bowtie \Hb \to \Hb\bowtie\Vc\,, \qquad
(g,h) \mapsto\Bigl(\bigl(h^{-1}\fiz\beta(g)^{-1}\bigr)^{-1},\bigl(\alpha(h)^{-1}\fde
g^{-1}\bigr)^{-1}\Bigr)\,;\\
\Hc\prode\Vb^{op}\to\Vb \prode \Hc^{op}\,, \qquad
 (x,\gamma)\mapsto (\gamma^{-1},x^{-1})\,.
\end{gather*}
\epf

We deduce that the constructions of doubles and duals commute:
\begin{corollary}\label{C:double-dual}
Let $(\Vc,\Hc)$ be a matched pair and consider its dual $(\Hc,\Vc)$.
We have
\[D(\Vc,\Hc)\cong D(\Hc,\Vc)\,.\]
\end{corollary}
\pf Take $\alpha$ and $\beta$ identities in Proposition~\ref{P:gendouble-dual}.
\epf

\section{Representations of matched pairs of groupoids}

\subsection{Definition and monoidal structure}

\

Let $\Pc$ be a set.  We denote by $\fpc$ the category of {\em quivers} over $\Pc$. The objects of
$\fpc$ are sets $X$ equipped with two maps $p:X\to\Pc$ and $q:X\to\Pc$, and the morphisms are maps $X\to
Y$ for which 
${\def\objectstyle{\scriptstyle}
\def\labelstyle{\scriptstyle}
\vcenter{\xymatrix@R-1.5pc@C-2pc{
X\ar[rd]_{p}\ar[rr] & & Y\ar[ld]^{p}\\ &{\Pc} }}}$ and
${\def\objectstyle{\scriptstyle}
\def\labelstyle{\scriptstyle}
\vcenter{\xymatrix@R-1.5pc@C-2pc{
X\ar[rd]_{q}\ar[rr] & & Y\ar[ld]^{q}\\ &{\Pc} }}}$
commute. This is a monoidal category: the tensor product of two objects $X$ and $Y$ 
is the set $X\Times{q}{p}Y$~\eqref{E:tensorcoprod} with $p(x,y)=p(x)$ and
$q(x,y)=q(y)$. The unit object is the set $\Pc$ with $p=q=\id:\Pc\to\Pc$.

\begin{definition}\label{D:representation}
Let $(\Vc,\Hc)$ be a matched pair of groupoids with base $\Pc$. A (set-theoretic) {\em representation} of
$(\Vc,\Hc)$ is an object $(\Ec,p,q)$ of $\fpc$ together with
\begin{itemize}
\item A left action  $\fde : \Hc \Times{r}{p} \Ec \to \Ec$ of $\Hc$ on $p$.
\item A  decomposition  $\Ec = \coprod_{g\in \Vc} \Ec(g)$ called the {\em grading} of $\Ec$ over $\Vc$. If
$e\in\Ec(g)$, we set $\abs{e}:= g$.
\end{itemize}
These are subject to the following conditions:
\begin{align}
\label{compcond-grading} p(e) = t(\abs{e}) & \text{ and } q(e)=b(\abs{e}) \text{ for all }e\in\Ec\,, \\
\label{compcond-action}
\abs{x\fde e} = x \fde \abs{e} & \text{ for all } (x, e) \in \Hc \Times{r}{p} \Ec\,.
\end{align}
\end{definition} 

Note that $p$ and $q$ are determined by the rest of the structure, in view of~\eqref{compcond-grading}.
The following diagram displays the above conditions:
\[\xymatrix@C+5pt{
{t(x\fde \abs{e})=l(x)}\ar@/_/[d]_{\abs{x\fde e}=x\fde \abs{e}}\ar@/^1pc/[r]^{x} &
{r(x)=p(e)=t(\abs{e})}\ar@/^/[d]^{\abs{e}}\\   {b(x\fde\abs{e})=l(x\fiz\abs{e})}\ar@/_1pc/[r]_{x\fiz
\abs{e}} & {r(x\fiz
\abs{e})=q(e)=b(\abs{e})} }\]

A morphism between representations $\Ec$ and $\F$ of $(\Vc,\Hc)$ is a map $\psi:\Ec\to\F$
which intertwines the actions of $\Hc$~\eqref{E:intertwiner} and preserves the gradings in the sense that
$\abs{\psi(e)}=\abs{e}$. It follows that $\psi$ preserves $p$ and $q$.

The resulting category of representations of the matched pair $(\Vc,\Hc)$ is denoted
$\Rep(\Vc,\Hc)$. It comes endowed with a functor $\ffl:\Rep(\Vc,\Hc)\to\fpc$
which forgets the action and the grading.

\begin{exa}\label{Ex:representations}  The
category of quivers over $\Pc$  may be identified with the category of representations of the initial
matched pair $(\Pc\times\Pc,\Pc)$  (Example~\ref{Ex:initial}). The action of $\Pc$ on a quiver
$(\Ec,p,q)$ is $P\fde e =e$, for $e\in \Ec$ with $p(e)=P$, and the grading over $\Pc\times\Pc$ is
$\abs{e}=\bigl(p(e),q(e)\bigr)$.

The category of representations of the terminal matched pair $(\Pc,\Pc\times\Pc)$
is less familiar. Its objects are sets $\Ec$ equipped with a 
decomposition $\Ec=\coprod_{P\in\Pc}\Ec_P$ and a family of maps $\pi_{P}:\Ec\to\Ec_P$, one for
each $P\in\Pc$, with the following properties: 
\[\pi_P(e)=e \text{  if $e\in\Ec_P$ and  }\pi_P\bigl(\pi_Q(e)\bigr)=\pi_P(e)\text{  for every
$e\in\Ec$, $P,Q\in\Pc$.}\]
In particular, each $\pi_P$ is a projection $\Ec\onto\Ec_P$. The structure of representation is as
follows. The quiver structure maps and the grading are $p(e)=q(e)=\abs{e}=P$, for $e\in\Ec_P$, while
the action of $\Pc\times\Pc$ is $(P,Q)\fde e=\pi_P(e)$, for $e\in\Ec_Q$.

More generally, a representation of a matched pair of the form $(\Vc,\Pc)$
(Example~\ref{Ex:semi-initial}) is just a set $\Ec$ with a grading over $\Vc$, and
a representation of a matched pair of the form $(\Pc,\Hc)$ is just a left action of
$\Hc$ on a map $p:\Ec\to\Pc$.

As a last example we describe the representations of the matched pair $M(X,Y)$ (Example~\ref{Ex:XY}).
These objects are quivers $\Ec$ with base $Y$, together with a decomposition 
$\Ec=\coprod_{x\in X}\Ec_x$ and a family of maps $\pi_{x}:\Ec\to\Ec_x$, one for
each $x\in X$, which are morphisms of quivers over $Y$ and have the following properties: 
\[\pi_x(e)=e \text{  if $e\in\Ec_x$ and  }\pi_x\bigl(\pi_{x'}(e)\bigr)=\pi_x(e)\text{  for every
$e\in\Ec$, $x,x'\in X$.}\]
Choosing $X$ or $Y$ as singletons we recover the above descriptions of the representations
of the initial and terminal matched pairs respectively.
\end{exa}

The notion of representation of $(\Vc,\Hc)$ makes use of the groupoid structure of $\Hc$ and
of the action of $\Hc$ on $\Vc$ only. 
Below, we introduce a monoidal structure on the category of representations which brings into play
the groupoid structure of $\Vc$ and the action of $\Vc$ on $\Hc$.

\begin{prop}\label{P:monope} Let  $\Ec$ and $\F$
 be representations of  a matched pair $(\Vc,\Hc)$. Consider the set
\[\Ec \otimes \F \ :=\   \Ec \Times{q}{p}  \F = \{(e,f) \in \Ec \times \F:
b(\abs{e}) = t(\abs{f})\}.\]
This set can be made into a representation of $(\Vc,\Hc)$ by defining
\begin{align}
\label{map-pt} p,q:\Ec \otimes \F\to\Pc\,, & \ \ p(e, f)=p(e)\,, \ \ q(e, f)=q(f)\,,\\
\label{action-pt} x \fde (e, f)  = & \bigl(x \fde e,  (x \fiz  \vert e \vert) \fde f\bigr)\,,\\
\label{grading-pt} \abs{(e,f)} = & \abs{e}\abs{f}\,.
\end{align}
\end{prop}
\begin{proof}
We check that this structure is well-defined and that $\Ec \otimes \F$
is indeed a representation of $(\Vc,\Hc)$.

Since $b(\vert e \vert)  = t(\vert f \vert)$, for $(e, f)\in \Ec \otimes \F$,
the grading is well-defined. Since $p(f) = b(\vert e \vert) =  r(x \fiz \vert e \vert)$,
the second component in \eqref{action-pt} is well-defined. Finally,
$b(\vert x\fde e \vert) = b( x  \fde \vert e \vert) = l( x \fiz \vert e \vert) = p((x \fiz  \vert e
\vert)\fde f)=t(\abs{(x \fiz  \vert e \vert)\fde f)})$; 
hence $(x \fde e,  (x \fiz  \vert e \vert) \fde f)\in \Ec \otimes \F$.

A straightforward computation shows that~\eqref{mp-0},~\eqref{mp-1}, and~\eqref{mp-1.5}
hold for \eqref{action-pt}.

We check condition~\eqref{compcond-grading} for $\Ec \otimes \F$: 
\[p(e,f)=p(e)=t(\abs{e})=t(\abs{e}\abs{f})=t\bigl(\abs{(e,f)}\bigr)\,,\]
by~\eqref{compcond-grading} for $\Ec$. Similarly, $q(e,f)=b\bigl(\abs{(e,f)}\bigr)$.
We check condition~\eqref{compcond-action} for $\Ec
\otimes
\F$: 
$$\abs{x \fde (e, f)}=\vert (x \fde e,  (x \fiz  \vert e \vert)
\fde f)\vert = \vert x \fde e \vert \vert (x \fiz  \vert e \vert) \fde f\vert
= (x \fde \vert e \vert)\bigl((x \fiz  \vert e \vert) \fde \vert f \vert\bigr)
=  x\fde (\vert e\vert \vert f \vert)\,,$$
by~\eqref{compcond-action} for $\F$ and by~\eqref{mp-3}. 
\end{proof}

The set $\Pc$ can be made into a representation of any matched pair $(\Vc,\Hc)$ of
groupoids with base $\Pc$. We set
$p=q=\id:\Pc\to\Pc$,  $\vert P\vert = \id_{\Vc} P\in \Vc$, and we let $\Hc$ act trivially on $p$, as
in~\eqref{E:trivialuno}, which in this case boils down to $x\fde r(x) = l(x)$ for all $x\in \Hc$. 

\begin{prop} \label{P:moncat} The operation $\Ec\otimes\F$ of Proposition~\ref{P:monope} turns
$\Rep(\Vc,\Hc)$ into a monoidal category. The unit object is $\Pc$, as above. 
The forgetful functor $\ffl:\Rep(\Vc,\Hc)\to\fpc$ is monoidal.
\end{prop}
\begin{proof}
If  $\Ec$, $\F$, and $\K$ are in $\Rep(\Vc,\Hc)$,
then there are canonical identifications 
$$(\Ec \otimes \F) \otimes \K \simeq \{(e, f, k)\in  \Ec \times \F \times \K:
q(e) = p(f)\,, \ q(f) = p(k)\}\simeq \Ec \otimes (\F \otimes \K)\,.$$
The grading is preserved under these identifications by associativity of the product of $\Vc$. As for the
action of $\Hc$, we have
\begin{multline*}
x\fde\bigl(e,(f,k)\bigr)=\bigl(x\fde e, (x\fiz\abs{e})\fde(f,k)\bigr)=
\Bigl(x\fde e, (x\fiz\abs{e})\fde f,\bigl((x\fiz\abs{e})\fiz{\abs{f}}\bigr)\fde k\Bigr)\\
=\Bigl(x\fde e,(x\fiz\abs{e})\fde f,\bigl(x\fiz(\abs{e}\,\abs{f})\bigr)\fde k\Bigr)=
\bigl(x\fde(e,f),(x\fiz\abs{(e,f)})\fde k\bigr)=x\fde\bigl((e,f),k\bigr)\,
\end{multline*}
as needed.
 
The verification of the remaining axioms offers no difficulty.
\end{proof}

\subsection{Some distinguished representations}\label{S:distinguished}

\

We discuss certain canonical representations of a matched pair $(\Vc,\Hc)$ which play an important
role in our proofs.

First, $\Vc$ becomes a representation of $(\Vc,\Hc)$ by means of the maps $t,b:\Vc\to\Pc$ (for $p$ and
$q$ respectively), the identity grading $\abs{g}=g$, and the action $\Hc\Times{r}{t}\Vc\map{\fde}\Vc$ from
the definition of matched pair. Conditions~\eqref{compcond-grading} and~\eqref{compcond-action} are
trivially satisfied.  This is in fact the final object in the category $\Rep(\Vc,\Hc)$, as one may easily
check:
\begin{lema}\label{L:terminal} Given any representation $\Ec \to \Pc$ of $(\Vc,\Hc)$,
there exists a unique morphism of representations $\psi: \Ec \to \Vc$. It is $\psi(e)=\abs{e}$.
\qed
\end{lema}

Another distinguished representation of  $(\Vc,\Hc)$
is $\Hc \Times{r}{t} \Vc$,  with action and grading given by
\begin{equation}\label{E:dist-rep}
x\fde (y, g) = (xy, g), \qquad \vert (x,g)\vert = x\fde g\,.
\end{equation}
Condition~\eqref{compcond-action} holds by~\eqref{mp-1}.
Condition~\eqref{compcond-grading} says that $p(x, g) = l(x)$ and $q(x,g)=l(x\fiz g)$.

This representation has several interesting subrepresentations.
 Let $\sigma: \Pc \to \Vc$ be a section of $t:\Vc\to\Pc$ and let
\[\Hc_{\sigma} := \{(x, \sigma r(x)): x\in \Hc\} \subseteq  \Hc \Times{r}{t} \Vc\,.\]
Then $\Hc_{\sigma}$ is a subrepresentation of $\Hc \Times{r}{t} \Vc$.

In particular, consider the section $\sigma_0$ defined by $\sigma_0(P) = \id_\Vc P$, $P\in \Pc$.
 We may identify $\Hc_{\sigma_0}$ with $\Hc$; then $p(x)=q(x)=l(x)$,
the action is $x\fde y=xy$, and the grading is $\vert x \vert = \id_\Vc l(x)$, by~\eqref{iduno}.

It is natural to wonder if there is a relation between the representations $\Vc$, $\Hc$, and $\Hc
\Times{r}{t}\Vc$ just defined.  One readily sees that $\Hc \Times{r}{t} \Vc$ is {\em not} the tensor
product of the representations $\Hc$ and $\Vc$. However, one has the following result.

\begin{lema}\label{morf}
The map $\morf: \Hc \Times{r}{t} \Vc \to \Vc \otimes \Hc
= \Vc \Times{b}{l} \Hc$ given by
\begin{equation}\label{isomorf}
\morf(x, g) = (x\fde g, x\fiz g)
\end{equation}
is an isomorphism of representations.
\end{lema}
\pf {}From~\eqref{E:cell} we see that $\morf$ is  well-defined and preserves
maps to $\Pc$: 
\[p(x,g)=l(x)=t(x\fde g)=p(x\fde g,x\fiz g)\,.\]
 The map $\varphi(g,x)=\bigl((x^{-1}\fiz g^{-1})^{-1},(g^{-1}\fde x^{-1})^{-1}\bigr)$ is the inverse of
$\morf$, by~\eqref{idcinco} and~\eqref{idseis}, so $\morf$ is a bijection.

We check that $\morf$ preserves actions. If  $x, y\in \Hc$ and $g\in \Vc$, with $r(x) = l(y)$, $r(y) =
t(g)$, then
\begin{multline*}
x\fde \morf(y, g) = x\fde (y \fde g, y\fiz g) \equal{\eqref{action-pt}}
(x\fde (y \fde g),(x\fiz (y \fde g))(y\fiz g))\\
 = (xy\fde g, xy\fiz g)=\morf(xy,g)\equal{\eqref{E:dist-rep}}\morf(x\fde (y, g)).
\end{multline*}
We check that $\morf$ preserves gradings. If  $x\in \Hc$ and
$g\in \Vc$, with $r(x) = t(g)$, then
\[\vert \morf(x, g)\vert = \vert (x \fde g, x\fiz g)\vert
\equal{\eqref{grading-pt}} (x \fde g) \id_{\Vc}l(x\fiz g) = (x \fde g) \equal{\eqref{E:dist-rep}} \vert
(x, g)\vert\,.\]
\epf

\bigbreak

There are additional important subrepresentations of $\Hc \Times{r}{t} \Vc$.
Recall that $\Hc_P = \{x\in \Hc: r(x) = P\}$, for $P\in \Pc$.
If $g\in \Vc$ is such that $t(g) = P$, then
$\Hc_P \times \{g\}$ is a subrepresentation of $\Hc \Times{r}{t} \Vc$. This representation
satisfies the following universal property.

\begin{lema}\label{expalema}
Let $\Ec$ be a representation of $(\Vc,\Hc)$. Let $e\in \Ec$, $g = \vert e \vert$, and $P =
t(g)$. There exists a unique morphism of representations $\rho_e: \Hc_P \times \{g\} \to \Ec$
 such that
\begin{equation}\label{expansion}
\rho_e(\id P, g) = e.
\end{equation}
\end{lema}

We refer to the morphism $\rho_e$ as the \emph{expansion} of $e$.

\pf  Define $\rho_e: \Hc_P \times \{g\} \to \Ec$ by $\rho_e(x, g) = x\fde e$.
This clearly satisfies \eqref{expansion}, and preserves action and grading.
Conversely, assume that $\rho: \Hc_P \times \{g\} \to \Ec$ is a morphism satisfying \eqref{expansion}.
Then $\rho(x, g) = \rho(x\id P, g) = x\fde \rho(\id P, g) = x\fde e = \rho_e(x, g)$, and uniqueness follows.
\epf

\bigbreak

Finally, we mention that $\Hc_P=\Hc_P \times \{\id P\}$
is a subrepresentation of $\Hc=\Hc_{\sigma_0}$, and
for any $y\in\Hc(P,Q)$, the right multiplication map $R_y: \Hc_P\to\Hc_Q$, $R_y(x) = xy$, is an
isomorphism of representations, with inverse $R_{y^{-1}}$.

\subsection{Restriction of representations}\label{S:restriction}

\

Let $(\alpha,\beta):(\Vb,\Hb)\to (\Vc,\Hc)$ be a morphism of matched pairs of groupoids with base $\Pc$.
We construct a functor 
\[\Res_\alpha^\beta:\Rep(\Vc,\Hc)\to\Rep(\Vb,\Hb)\]
called the {\em restriction along} $(\alpha,\beta)$. 
Let $\Ec$ be a representation  of $(\Vc,\Hc)$ with structure maps $p,q:\Ec \to \Pc$, action $\fde$ and
grading $\abs{\ }$. Consider the  action of $\Hb$ on the map $p:\Ec \to\Pc$ given by
\begin{equation}\label{functor3}
h \fde e : = \alpha(h) \fde e\,, \text{ for all $(h,e)\in \Hb\Times{r}{p} \Ec$,} 
\end{equation}
and the grading $\Vert \ \Vert: \Ec \to \Vb$  given by 
\begin{equation}\label{functor4}
\Vert e \Vert := \beta(\vert e \vert)\,.
\end{equation}
We let $\Res^\beta_\alpha(\Ec)$ denote
the set $\Ec$, equipped with the same maps $p,q:\Ec \to \Pc$ and these new action and grading.

Given a morphism of representations $\psi:\Ec\to\F$, we let $\Res^\beta_\alpha(\Ec)(\psi)=\psi$.

\begin{prop}\label{restr} Let $(\alpha,\beta):(\Vb,\Hb)\to (\Vc,\Hc)$ be a morphism of matched pairs of 
groupoids.
\begin{itemize}
\item[(a)] For any representation $\Ec$ of $(\Vc,\Hc)$,
$\Res^\beta_\alpha(\Ec)$ is a representation of $(\Vb,\Hb)$.
\item[(b)] $\Res^\beta_\alpha:\Rep(\Vc,\Hc)\to\Rep(\Vb,\Hb)$ is a monoidal functor and diagram
\[\xymatrix@R-1pc@C-1.5pc{
{\Rep(\Vc,\Hc)}\ar[rd]_{\ffl}\ar[rr]^{\Res^\beta_\alpha} & & {\Rep(\Vb,\Hb)}\ar[ld]^{\ffl}\\ &{\fpc} }\]
commutes.
\item[(c)] If $(\chi, \omega):(\Vc,\Hc)\to(\Vg,\Hg)$ is another morphism of matched pairs then diagram
\[\xymatrix@R-1pc@C-1.5pc{
{\Rep(\Vg,\Hg)}\ar[rd]_{\Res^\omega_\chi}\ar[rr]^{\Res^{\beta\circ\omega}_{\chi\circ\alpha}} & &
{\Rep(\Vb,\Hb)}\ar[ld]^{\Res^\beta_\alpha}\\ &{\Rep(\Vc,\Hc)}\ar[ur]_{\Res^\beta_\alpha} }\] 
commutes.
\end{itemize}
\end{prop}
\pf (a) We first check that the new action and grading satisfy~\eqref{compcond-grading}
and~\eqref{compcond-action}. Given $e\in\Ec$ we have
\[p(e)=t(\abs{e})=t\bigl(\beta(\abs{e})\bigr)=t(\norm{e})\,,\]
by~\eqref{compcond-grading} for $\Ec$ as a representation of $(\Vc,\Hc)$.
Similarly, $q(e)=b(\norm{e})$. If $(h,e)\in\Hb\Times{r}{p} \Ec$ then
$$
\Vert h\fde e \Vert = \Vert \alpha (h) \fde e \Vert = \beta(\vert \alpha (h) \fde e \vert)
= \beta(\alpha (h) \fde \vert  e \vert) \equal{\eqref{functor1}}  h \fde \beta(\vert  e \vert) = h\fde
\Vert  e\Vert\,;$$
by~\eqref{compcond-action} for $\Ec$ as a representation of $(\Vc,\Hc)$.  Thus, 
$\Res^\beta_\alpha$ is a well-defined functor, the rest of the verifications being straightfoward.

We next check that $\Res^\beta_\alpha$ is monoidal. Since the functor does not change the
maps $p$ and $q$, $\Res^\beta_\alpha(\Ec\otimes\F) =\Res^\beta_\alpha(\Ec) \otimes \Res^\beta_\alpha(\F)$
as sets. Let $e\in \Ec$, $f\in \F$. The degree of $(e,f)$ in $\Res^\beta_\alpha(\Ec \otimes \F)$ is
$$
\Vert (e, f) \Vert =  \beta(\vert (e, f)  \vert) =   \beta(\vert e  \vert \vert f  \vert)
=   \beta(\vert e  \vert)  \beta(\vert f  \vert) = \Vert e  \Vert \Vert f  \Vert,
$$
which agrees with the degree of $(e,f)$ in $\Res^\beta_\alpha(\Ec) \otimes \Res^\beta_\alpha(\F)$.
 Let  $h\in \Hb$. The action of $h$ on $(e,f)$ in $\Res^\beta_\alpha(\Ec \otimes \F)$ is
\begin{multline*}
h\fde (e,f) = \alpha(h) \fde (e,f) =  \Bigl(\alpha(h) \fde e, \bigl(\alpha(h) \fiz \vert e \vert\bigr)
\fde f\Bigr)
\\ \equal{\eqref{functor2}} \Bigl(\alpha(h) \fde e, \alpha\bigl(h \fiz \beta(\vert e \vert)\bigl) \fde
f\Bigr) = \Bigl(\alpha(h) \fde e, \alpha\bigl(h \fiz \Vert e \Vert\bigr) \fde f\Bigr),
\end{multline*}
which agrees with the action of $h$ on $(e,f)$ in $\Res^\beta_\alpha(\Ec) \otimes \Res^\beta_\alpha(\F)$.

The proof of (c) is left to the reader.
\epf

\begin{theorem}\label{restr-uniq}
 Let $F:\Rep(\Vc,\Hc)\to\Rep(\Vb,\Hb)$ be a monoidal functor such that diagram
\[\xymatrix@R-1.5pc@C-2pc{
{\Rep(\Vc,\Hc)}\ar[rd]_{\ffl}\ar[rr]^{F} & & {\Rep(\Vb,\Hb)}\ar[ld]^{\ffl}\\ &{\fpc} }\]
commutes. There exists a unique morphism of matched pairs of groupoids $(\alpha,\beta):(\Vb,\Hb)\to
(\Vc,\Hc)$ such that
\[F = \Res^\beta_\alpha\,.\]
\end{theorem}
\pf We make use of the representation $\Vc$ of $(\Vc,\Hc)$ of Lemma~\ref{L:terminal}.
Since $F(\Vc)$ is an object in $\Rep(\Vb,\Hb)$
with underlying set $\Vc$, we may define a map $\beta:\Vc\to\Vb$ by
\[\beta(v) = \Vert v \Vert\,.\] 
Let $\Ec$ be an arbitrary representation
of $(\Vc,\Hc)$. By Lemma~\ref{L:terminal}, the map $\psi:\Ec\to \Vc$ defined by $\psi(e)=\abs{e}$
 is a morphism in $\Rep(\Vc,\Hc)$. Since $F$ preserves forgetful functors, the same map $\psi$ is a
morphism $F(\Ec)\to F(\Vc)$ in $\Rep(\Vb,\Hb)$; in particular, it must preserve the gradings. We
deduce that
\begin{equation*}
\beta(\abs{e})=\beta(\psi(e))=\norm{\psi(e)}=\norm{e}\,.
\end{equation*}
We have obtained condition~\eqref{functor4} in the definition of restriction functor.

{}From $F(\Vc \otimes \Vc) =  F(\Vc) \otimes F(\Vc)$ we deduce that
$\beta$ is a morphism of groupoids.

Let $h\in \Hb$, $P = r(h)$, and $g\in \Vc$ be such that $t(g) = P$.
Consider the representation $\Hc_P \times \{g\}$ of $(\Vc,\Hc)$.
 Using the action of $\Hb$ on $F(\Hc_P \times \{g\})$, we define
an element $\alpha_g(h) \in \Hc_P$ by
$$
h\fde (\id_{\Hc} P, g) = (\alpha_g(h), g).
$$
Let  $\Ec$ be a representation of $(\Vc,\Hc)$. Let $e\in \Ec$ be such that
$g = \vert e \vert$. We claim that the action of $\Hb$ on $F(\Ec)$ is given by
$h\fde e = \alpha_g(h) \fde e$. To see this, consider the expansion $\rho_e$ of $e$, that is the map
$\rho_e:\Hc_P \times \{g\}\to\Ec$, $\rho_e(x)=x\fde e$. By Lemma~\ref{expalema}, $\rho_e$ is a
morphism in $\Rep(\Vc,\Hc)$. Since $F$ preserves forgetful functors,
$\rho_e$ is also a morphism in $\Rep(\Vb,\Hb)$; in particular, it preserves the actions of $\Hb$.
Therefore,
\begin{equation}\tag{$*$}
h\fde e = h\fde \bigl(\rho_e(\id P, g)\bigr) = \rho_e\bigl(h\fde (\id P, g)\bigr) = \rho_e(\alpha_g(h), g)
= \alpha_g(h) \fde e\,.
\end{equation}

We next show that $\alpha_g(h)$ is independent of $g$. Let $f\in\Vc$ be another element such
that $t(f)=P$. Let $X = \Hc_P \times \{f\}$
and $Y = \Hc_P \times \{g\}$. We compute
the action of $\Hb$ on $F(X \otimes Y) = F(X) \otimes F(Y)$ in two different ways.
Let $\bigl((x, f),(y,g)\bigr)\in X\otimes Y$. According to~\eqref{grading-pt} and~\eqref{E:dist-rep}, the
degree of this element is
\[\abs{((x, f),(y,g)\bigr)}=\abs{(x,f)}\abs{(y,g)}=(x\fde f)(y\fde g)\,.\]
The top of $(x\fde f)(y\fde g)\in\Vc$ is, by definition of $\Hc_P$, $t(f)=P$. Therefore, we may apply
$(*)$ to this element and conclude that the action of $h$ on $F(X \otimes Y)$ is
\begin{align*}
h\fde \bigl((x, f),(y, g)\bigr)
&= \alpha_{(x\fde f)(y\fde g)}(h) \fde \bigl((x, f),(y, g)\bigr)\\
&=  \Biggl(\Bigl(\alpha_{(x\fde f)(y\fde g)}(h)x,  f\Bigr),\Bigl(\bigl(\alpha_{(x\fde f)(y\fde g)}(h) \fiz
(x\fde f)\bigr) y   , g\Bigr)\Biggr)\,.
\end{align*}
On the other hand, the action of $h$ on $F(X) \otimes F(Y)$ is
\begin{align*}
h\fde \bigl((x, f),(y, g)\bigr)
&=\Bigl(h\fde (x, f),\bigl(h\fiz \Vert (x, f)\Vert\bigr)\fde (y, g)\Bigr)\\
&= \Biggl(\bigl(\alpha_{x\fde f}(h) x,  f\bigr),\Bigl(\alpha_{y\fde g}\bigl(h\fiz \beta(x \fde
f)\bigr) y , g\Bigr)\Biggr)\,.
\end{align*}
Comparing these two expressions we obtain
\begin{align} \label{comp1}
\alpha_{(x\fde f)(y\fde g)}(h)&=\alpha_{x\fde f}(h)\,,
\\ \label{comp2}
\alpha_{(x\fde f)(y\fde g)}(h) \fiz (x\fde f) &=  \alpha_{y\fde g}\bigl(h\fiz \beta(x \fde f)\bigr)\,.
\end{align}

Choose $x=y = \id_{\Hc} P\in\Hc_P$ and $f = \id_{\Vc} P$. Since
$q((x,f)) = P = p((y,g))$, we have $\bigl((x,f),(y,g)\bigr)\in X\otimes Y$.
Hence,~\eqref{comp1} applies and we obtain that $\alpha_{g}(h) = \alpha_{\id_{\Vc} P}(h)$. Thus,  we may
define a map
$\alpha: \Hb \to \Hc$ by $\alpha(h)=\alpha_g(h)$ where $g$ is any element of $\Vc$ such that $t(g)=r(h)$.
Equation $(*)$ becomes condition~\eqref{functor3} in the definition of restriction functor.

It follows  easily now that $\alpha$ is a morphism of groupoids:
if $h,k\in \Hb$, $r(k) = l(h)$ then $(*)$ implies 
$$
(\alpha(kh), \id_{\Vc}  P) = kh \fde (\id_{\Hc} , \id_{\Vc}  P) =
k\fde\left(h \fde (\id_{\Hc} , \id_{\Vc}  P)\right) =
k\fde(\alpha(h) , \id_{\Vc}  P) = (\alpha(k)\alpha(h) , \id_{\Vc}  P).
$$
Also,~\eqref{comp2} yields $\alpha (h) \fiz (x\fde f) =  \alpha\bigl(h\fiz \beta(x \fde f)\bigr)$.
Choosing $x=\id_{\Hc}(t(f))$ we deduce~\eqref{functor2} in the definition of morphism of matched
pairs of groupoids. 
To complete the proof, it only
remains to check that~\eqref{functor1} holds. Let $h\in \Hb$, $g\in \Vc$. Then
$$
\beta (\alpha (h) \fde g) = \Vert \alpha (h) \fde g\Vert= \Vert h \fde g\Vert =  h \fde \Vert g\Vert
= h \fde \beta(g)\,.
$$
Thus, $(\alpha,\beta):(\Vb,\Hb)\to(\Vc,\Hc)$ is a morphism of matched pairs and 
$\Res^\beta_\alpha=F$.

Uniqueness can be traced back through the above argument.
\epf

\begin{exa}\label{Ex:forgetful} Consider the forgetful functor $\ffl:\Rep(\Vc,\Hc)\to\fpc$. 
As discussed in Example~\ref{Ex:representations}, $\fpc=\Rep(\Pc\times\Pc,\Pc)$.
 According to Theorem~\ref{restr-uniq}, $\ffl$ must be the
 restriction functor along a morphism
of matched pairs $(\Pc\times\Pc,\Pc)\to (\Vc,\Hc)$. In fact, one has 
$\ffl=\Res_{\alpha_0}^{\beta_0}$ where $(\alpha_0,\beta_0)$ is the morphism of 
Proposition~\ref{P:initial}.

Dually, consider the  unique morphism of matched pairs from $(\Vc,\Hc)$ to
$(\Pc,\Pc\times\Pc)$. The corresponding restriction functor allows to view any
representation of $(\Pc,\Pc\times\Pc)$ (see Example~\ref{Ex:representations}) as a representation of
$(\Vc,\Hc)$. Such representations are called {\em trivial}. None of the representations of
Section~\ref{S:distinguished} are trivial.
\end{exa}

\section{Centers and centralizers}\label{S:centers}

\subsection{Generalities on centralizers, centers, and braidings}\label{S:generalities}

\ 

We review the connection between centers and braidings for monoidal categories, and discuss the
related notion of centralizers. We restrict our considerations to the case of strict monoidal functors
between strict monoidal categories, as this is the case of present interest.
For more information,  see~\cite{js} or~\cite[XIII.4]{kas}.

Let $F: \C \to \D$ be a (strict) monoidal functor
between (strict) monoidal categories. There is a
(strict) monoidal category $Z(F)$, the {\em centralizer} of $F$, defined as follows.

\begin{itemize}
\item The objects of $Z(F)$ are pairs $(D, \zeta)$, where $D$ is an object of $\D$
and $\zeta_{(-)}:  F(-) \otimes D\to D \otimes F(-) $
is a natural isomorphism such that the following diagrams commute:
\begin{equation}\label{cond1centr}
\xymatrix@C-30pt@R-7pt{
{F(X\otimes Y)\otimes D}\ar[rr]^{\zeta_{X\otimes Y}}\ar@{=}[d] & & {D\otimes F(X\otimes Y)}\ar@{=}[d]\\
{F(X)\otimes F(Y)\otimes D}\ar[rd]_{\id\otimes \zeta_Y} & & {D\otimes F(X)\otimes F(Y)}\\
& {F(X)\otimes D\otimes F(Y)}\ar[ur]_{\zeta_X\otimes\id} 
}
\end{equation}
\begin{equation}\label{cond2centr}
\xymatrix@C-10pt@R-10pt{
{F(\uno)\otimes D}\ar[rr]^{\zeta_{\uno}}\ar@{=}[d] & & {D\otimes F(\uno)}\ar@{=}[d] \\
{\uno\otimes D}\ar@{=}[r] & D\ar@{=}[r] & {D\otimes \uno}
}
\end{equation}

\item A morphism $f: (D, \zeta)\to (E, \vartheta)$ in $Z(F)$ is a map
$f: D \to E$ in $\D$ such that for all objects $X$ of $\C$ the following diagram commutes 
\[\xymatrix@C-5pt@R-10pt{
{F(X)\otimes D}\ar[r]^{\zeta_X}\ar[d]_{\id\otimes f} & {D\otimes F(X)}\ar[d]^{f\otimes\id}\\
{F(X)\otimes E}\ar[r]_{\vartheta_X} & {E\otimes F(X)}
}\] 

\item The tensor product of two objects $(D, \zeta)$ and $(E, \vartheta)$ in $Z(F)$ is the object
$(D\otimes E, \zeta\otimes\vartheta)$ where $(\zeta\otimes\vartheta)_{(-)}:D\otimes E\otimes F(-)\to
F(-)\otimes D\otimes E$ is defined for each object $X$ of $\C$ by the diagram 
\begin{equation}\label{ptcentr}
\xymatrix@C-20pt@R-10pt{
{D\otimes E\otimes F(X)}\ar[rd]_{\id\otimes \vartheta_X} \ar@{-->}[rr]^{(\zeta\otimes\vartheta)_{X}}
& & {F(X)\otimes D\otimes E}\\ & {D\otimes F(X)\otimes E}\ar[ru]_{\zeta_X\otimes\id} & 
}
\end{equation}

\item The tensor product of morphisms is the tensor product of morphisms in $\D$.

\item The unit object of $Z(F)$ is $(\uno,\id)$ where  
$\id_{(-)}:\uno\times F(-)\to F(-)\times \uno$
is the identity natural transformation.

\end{itemize}

There is a strict monoidal functor $U:Z(F)\to\D$ given by $U(D,\zeta)=D$.

\bigbreak
The {\em center} of $\C$ is the centralizer of the identity functor $\C\to\C$.
It is denoted $Z(\C)$. It is a \emph{braided} category with braiding
\begin{equation}\label{E:braiding-center}
c_{(E, \vartheta),(D, \zeta)}: (E, \vartheta)  \otimes (D, \zeta) \rightarrow (D, \zeta) \otimes (E, \vartheta)\,, 
\qquad c_{(E, \vartheta),(D, \zeta)} = \zeta_E\,.
\end{equation}

Let $U: Z(\C) \to \C$ be as before. 
Braidings in $\C$ are in one-to-one correspondance with monoidal sections of $U$. More precisely,
a family of maps $c_{E,D}:E\otimes D\to D\otimes E$ is a braiding in the monoidal category $\C$ if and
only if  $D\mapsto (D,c_{-,D})$ defines a monoidal functor $\C\to Z(\C)$.

\subsection{Centralizers of restriction functors}\label{S:centralizers}

\

Let $(\alpha, \beta):(\Vb,\Hb)\to(\Vc,\Hc)$ be a morphism of matched
pairs of groupoids. Consider the restriction functor $\Res^\beta_\alpha:\Rep(\Vc,\Hc)\to\Rep(\Vb,\Hb)$
(Section~\ref{S:restriction}). We show that the centralizer of this monoidal functor is the category
of representations of a matched pair; namely, the generalized double of $(\alpha, \beta)$ (Definition
~\ref{D:double}).

\begin{theorem}\label{T:centralizer} Let $(\alpha, \beta):(\Vb,\Hb)\to(\Vc,\Hc)$ be a morphism of matched
pairs of groupoids with base $\Pc$. There is a strict monoidal equivalence
\[ \Rep D(\alpha,\beta)\cong Z(\Res^\beta_\alpha)\,,\]
and the following diagram is commutative
\[\xymatrix{
{\Rep D(\alpha,\beta)}\ar[r]^{\cong}\ar[rd]_{\ffl}& {Z(\Res^\beta_\alpha)}\ar[r]^{U}&
{\Rep(\Vb,\Hb)}\ar[ld]^{\ffl}\\ 
&{\Quiv(\Pc)}& \\
}\]
\end{theorem}

We approach the proof of this theorem in several steps.

Let $(D, \zeta)$ be an object in $Z(\Res^\beta_\alpha)$. In particular,
$D$ is a representation of $(\Vb,\Hb)$ and a quiver over $\Pc$. 
Let $p,q:D\to\Pc$ be the quiver structure maps and $\norm{d}\in\Vb$ denote the degree of an element
$d\in D$. 

For any object $\Ec$ in $\Rep(\Vc,\Hc)$, we have the natural isomorphism 
$\zeta_{\Ec}:\Res^\beta_\alpha(\Ec)\otimes D\to D\otimes\Res^\beta_\alpha(\Ec)$.

For any sets $X$ and $Y$, we let $\pr_1$ and $\pr_2$ denote the  projections from $X\times Y$ onto $X$
and $Y$, respectively.

\begin{claim} There is a map  $\rightharpoonup: \Vc\Times{b}{p} D \to D$ such that
\begin{equation}\label{action12}
\pr_1\zeta_{\Ec} (e, d)  =  \vert e\vert \rightharpoonup d\,, 
\end{equation}
for all  $e\in \Ec$ and $d\in D$ with $b(\vert e \vert) =p(d)$.
Moreover, $\rightharpoonup$ is a left action of $\Vc$ on  
$p:D\to\Pc$.
\end{claim}
\pf Consider the object $\Vc$ of $\Rep(\Vc,\Hc)$ of Lemma~\ref{L:terminal}. Define
 $$g \rightharpoonup d := \pr_1 \zeta_{\Vc} (g, d)\in D\,,$$
for $g\in \Vc$ and $d\in D$ with $b(g) =p(d)$.
Applying the naturality condition of $\zeta$ to the unique morphism of representations $\Ec\to \Vc$
(the degree) we deduce~\eqref{action12}.

We verify the conditions for $\rightharpoonup: \Vc \Times{b}{p} D \to D$ to be 
a left action of $\Vc$ on $p:D\to\Pc$: \eqref{mp-0} holds since $\zeta_{\Vc}$ commutes with maps
to $\Pc$ (being a map of quivers over $\Pc$), \eqref{mp-1} follows from~\eqref{cond1centr} applied to
$X=Y=\Vc$, and~\eqref{mp-1.5} follows from~\eqref{cond2centr}.
\epf
\begin{claim} There is a map  $[\ ]:D\to\Hc$ such that 
\begin{equation}\label{action13}
\zeta_{\Ec} (e, d) =\bigl(\abs{e} \rightharpoonup d, ([d] \fiz \abs{e}^{-1}) \fde e\bigr)\,,
\end{equation}
for all  $e\in \Ec$ and $d\in D$ with $b(\vert e \vert) =p(d)$.
Moreover, the map $d\mapsto (\norm{d},[d])$ maps $D$ to $\Vb \prode \Hc^{op}$.
\end{claim}
\pf Consider the object $\Hc$ of $\Rep(\Vc,\Hc)$ of Section~\ref{S:distinguished}. Define,
for $d\in D$,
\[[d] := \pr_2 \zeta_{\Hc} \bigl(\id_{\Hc} p(d), d\bigr)\in \Hc\,.\]
Since $\zeta_{\Hc} \bigl(\id_{\Hc} p(d), d\bigr)\in D\otimes\Res^\beta_\alpha(\Hc) =D\Times{q}{l}\Hc$, we
have
\begin{equation*}
l([d]) = q\Bigl(\pr_1 \zeta_{\Hc} \bigl(\id_{\Hc} p(d), d\bigr)\Bigr)
\equal{\eqref{action12}}
q\bigl(\abs{\id_{\Hc} p(d)} \rightharpoonup d\bigr) =q\bigl(\id_{\Vc} p(d) \rightharpoonup d\bigr)
\equal{\eqref{mp-1.5}}q(d)\equal{\eqref{compcond-grading}} b(\Vert d\Vert)\,.
\end{equation*}
(we used~\eqref{compcond-grading} for $D$ as a representation of $(\Vb,\Hb)$). 

Let $P:=p(d)\in\Pc$. Consider the subobject $\Hc_P$ of $\Hc$ in $\Rep(\Vc,\Hc)$
(Section~\ref{S:distinguished}). By naturality of $\zeta$,
\[ [d]=\pr_2 \zeta_{\Hc_P} \bigl(\id_{\Hc} P, d\bigr)\in \Hc_P\,.\]
Hence, by definition of $\Hc_P$, 
\[r([d])=P=p(d) \equal{\eqref{compcond-grading}} t(\Vert d\Vert)\,.\]
Thus $(\norm{d},[d])$ belongs to $\Vb\prode\Hc^{op}=(\Vb\Times{b}{l}\Hc)\cap(\Vb\Times{t}{r}\Hc)$.

Let $e$ be as above, $g=\abs{e}\in\Vc$, and $Q=t(g)$. Consider the object
$\Hc_{Q}\times\{g\}$ in $\Rep(\Vc,\Hc)$ (Section~\ref{S:distinguished}). 
Since $\abs{(\id_{\Hc} Q, g)}=g$,
~\eqref{action12} says that
\[\zeta_{\Hc_{Q}\times \{g\}} (\id_{\Hc} Q, g, d) = (g \rightharpoonup d, [d]_g, g)\]
for some element $[d]_g \in \Hc_Q$.

Let $\rho_e: \Hc_Q \times \{g\} \to \Ec$ be the  expansion of $e$ (Lemma
\ref{expalema}). Applying the naturality of $\zeta$ to the morphism $\rho_e$ we deduce that
 \begin{equation}\tag{$\ast$}
 \zeta_{\Ec} (e, d)  = (g \rightharpoonup d,[d]_g \fde e)\,.
\end{equation}
Consider now the isomorphism $\morf$ of Lemma \ref{morf}. Combining the naturality of $\zeta$
with~\eqref{cond1centr} we obtain the commutativity of
\[\xymatrix@C-20pt{ 
{\Res^\beta_\alpha(\Hc \Times{r}{t} \Vc)  \otimes  D}
\ar[rr]^{\zeta_{\Hc \Times{r}{t}\Vc}}\ar[d]_{\Res^\beta_\alpha(\morf)  \otimes  \id} & & {D \otimes
\Res^\beta_\alpha(\Hc \Times{r}{t} \Vc)}\ar[d]^{\id\otimes\Res^\beta_\alpha(\morf)} \\
{\Res^\beta_\alpha(\Vc) \otimes \Res^\beta_\alpha(\Hc)  \otimes D}
\ar[rd]_{\id \otimes \zeta_{\Hc}} & &  {D \otimes  \Res^\beta_\alpha(\Vc)
\otimes\Res^\beta_\alpha(\Hc)}\\
& {\Res^\beta_\alpha(\Vc) \otimes D\otimes \Res^\beta_\alpha(\Hc)}
\ar[ru]_{\zeta_{\Vc} \otimes \id} & \\
}\]
On one hand,
$$
\left(\id \otimes \Res^\beta_\alpha(\morf) \right) \zeta_{\Hc \Times{r}{t} \Vc} (\id_{\Hc} Q, g, d)
\equal{\eqref{action13}} \left(\id \otimes \Res^\beta_\alpha(\morf)\right) (g \rightharpoonup d, [d]_g,g)
\equal{\eqref{isomorf}} (g \rightharpoonup d, [d]_g \fde g, [d]_g \fiz g)\,.$$
On the other, by~\eqref{isomorf} and~\eqref{iddos},
\begin{multline*}
(\zeta_{\Vc} \otimes \id)(\id \otimes \zeta_{\Hc}) (\Res^\beta_\alpha(\morf)\otimes \id ) 
(\id_{\Hc} Q, g, d)= (\zeta_{\Vc} \otimes \id)(\id \otimes \zeta_{\Hc}) (g, \id_{\Hc} P, d)\\ 
=  (\zeta_{\Vc} \otimes \id) (g, d, [d])
\equal{(\ast)}  (g \rightharpoonup d, [d]_g \fde g, [d]).
\end{multline*}
Hence $[d] = [d]_g \fiz g$, or $[d]_g = [d] \fiz g^{-1}$. Together with $(\ast)$,
this proves~\eqref{action13}. 
\epf

\begin{claim} The action $\rightharpoonup$ of $\Vc$ on $D$ and the action of $\Hb$ on $D$ extend to
an action of $\Vc \bowtie \Hb$ on $D$.  Together with the map $(\norm{\ },[\ ])$, they turn
$D$ into a representation of the matched pair $D(\alpha,\beta)=(\Vb \prode \Hc^{op},\Vc \bowtie \Hb)$.
\end{claim}
\pf We make use of the fact that $\zeta_{\Ec}$ preserves the action of $\Hb$. Let $h\in\Hb$. On one hand,
\begin{multline*}
\zeta_{\Ec}\bigl(h\fde (e,d)\bigr) \equal{\eqref{action-pt}}
\zeta_{\Ec}\bigl(h\fde e , (h\fiz \Vert e \Vert) \fde d\bigr)
 \equal{\eqref{functor3},\eqref{functor4}} \zeta_{\Ec}\Bigl(\alpha(h)\fde e, \bigl(h\fiz \beta(
\vert e\vert)\bigr)\fde d\Bigr)\\ 
\equal{\eqref{action13}} \Biggl(\vert\alpha(h)\fde e\vert \rightharpoonup 
\Bigl(\bigl(h\fiz \beta( \vert e \vert)\bigr) \fde d\Bigr), \Bigl(\Bigl[\bigl(h\fiz \beta( \vert e
\vert)\bigr)
\fde d\Bigr]
\fiz\vert\alpha(h)\fde e\vert^{-1}\Bigr) \fde \bigl(\alpha(h)\fde e\bigr)\Biggr)\,.
\end{multline*} 
On the other,
\begin{multline*}
h\fde \zeta_{\Ec} (e,d)  \equal{\eqref{action13}} 
h\fde \bigl(\vert e\vert \rightharpoonup d, ([d] \fiz\vert e\vert^{-1})\fde e\bigr)
\equal{\eqref{action-pt}} \Bigl(h\fde (\vert e \vert \rightharpoonup d),
\bigl(h \fiz \Norm{\,\vert e \vert \rightharpoonup d\,}\bigr) \fde \bigl([d] \fiz \vert e \vert^{-1}\bigr)
\fde e\Bigr)\\
\equal{\eqref{functor3}}\Bigl(h\fde (\vert e \vert \rightharpoonup d),
\alpha\bigl(h \fiz \Norm{\,\vert e \vert \rightharpoonup d\,}\bigr) \fde \bigl([d] \fiz \vert e
\vert^{-1}\bigr)
\fde e\Bigr) \,.
\end{multline*}
The equality between the first components of the above expressions and
formula~\eqref{E:diag-product-2} for the product of the groupoid $\Vc \bowtie \Hb$ imply that
the actions $\rightharpoonup$ and $\fde$ combine to give a left action of $\Vc \bowtie \Hb$
 on  $D$.

Consider the equality between the second components, in the special case when $\Ec=\Hc$ and $e$ is an
identity. Recalling the structure of $\Hc$ as a representation of $(\Vc,\Hc)$
(Section~\ref{S:distinguished}), we obtain
\begin{equation}\label{action14}
[h\fde d]\alpha(h)=
\alpha\bigl(h \fiz \norm{d}\bigr)[d]\,.
\end{equation}

We need to derive additional identities. Applying~\eqref{cond1centr} to the product
$\Res^\beta_\alpha(\Hc)\otimes\Res^\beta_\alpha(\Vc)\otimes D$ and the element
$(\id_{\Hc} Q, g,d)$ we obtain
\begin{equation}\label{action14.5}
[g\rightharpoonup d]=[d]\fiz g^{-1}\,.
\end{equation}

Finally, we make use of the fact that $\zeta_{\Ec}$ preserves degrees (in $\Vb$). Calculating degrees
on both sides of~\eqref{action13} (and letting $g=\abs{e}$) we obtain
\begin{equation}\label{action15}
\beta(g)\Vert d \Vert = \Vert g\rightharpoonup d\Vert\beta\bigl(([d] \fiz g^{-1}) \fde g\bigr)\,, 
\end{equation}
for any $d\in $ and $g \in\Vc$ with $b(g)=q(d)$.

We have already seen that $D$, with its given structure of quiver over $\Pc$, carries a grading
over $\Vb\prode\Hc^{op}$ (Claim 2) and a left action of $\Vc \bowtie \Hb$. To conclude that
$D$ is a representation of the matched pair $D(\alpha,\beta)=(\Vb\prode\Hc^{op},\Vc \bowtie \Hb)$
we only need to check conditions~\eqref{compcond-grading} and~\eqref{compcond-action}.
The former boils down to the corresponding condition for $D$ as a representation of
$(\Vb,\Hb)$. The latter  requires the verification of the following two identities (in view
of~\eqref{action11}):
\begin{align}
\bigl[g\rightharpoonup(h\fde d)\bigr] &=
\Bigl(\alpha\bigl(h\fiz\norm{d}\bigr)[d]\alpha(h)^{-1}\Bigr)\fiz g^{-1}\,, \label{action16}\\
\norm{g\rightharpoonup(h\fde d)} &=
\beta(g)\bigl(h\fiz\norm{d}\bigr)\beta\Bigl(\bigl(\alpha(h\fiz\norm{d})[d]\alpha(h)^{-1}\bigr)\fiz
g^{-1}\Bigr)\,.\label{action17}
\end{align}
Now, \eqref{action16} is a consequence of~\eqref{action14} and~\eqref{action14.5}, while~\eqref{action17}
follows from~\eqref{idtres}, ~\eqref{action15}, ~\eqref{action16}, and~\eqref{compcond-action} for $D$ as
a representation of $(\Vb,\Hb)$.
\epf

We have shown that any object $D$ in the centralizer $Z(\Res^\beta_\alpha)$ can be endowed with the
structure of a representation of the matched pair $D(\alpha,\beta)$. Conversely, any such representation
$D$ gives rise to an object in $Z(\Res^\beta_\alpha)$. First, $D$ is a representation of
$(\Vb,\Hb)$ by restriction via the morphism of matched pairs $(\Vb,\Hb)\to D(\alpha,\beta)$ of
Remark~\ref{R:inclusion}. Then, one defines a transformation $\zeta_{\Ec}:\Res^\beta_\alpha(\Ec)\otimes
D\to D\otimes\Res^\beta_\alpha(\Ec)$ by means of~\eqref{action13}. Reversing some of the arguments
in the above proofs, one verifies that $(D,\zeta)$ is an object in $Z(\Res^\beta_\alpha)$.

It is clear that these constructions yield an equivalence of categories $ \Rep
D(\alpha,\beta)\cong Z(\Res^\beta_\alpha)$, whose composition with the functor
$U:Z(\Res^\beta_\alpha)\to\Rep(\Vb,\Hb)$ preserves the forgetful functors to $\Quiv(\Pc)$.
 To complete the proof of Theorem~\ref{T:centralizer}, we need to verify that
the monoidal structures are preserved under this equivalence.

\begin{claim}  Tensor products of objects in the categories $Z(\Res^\beta_\alpha)$ and $\Rep
D(\alpha,\beta)$ agree. 
\end{claim}
\pf Let $D$ and $D'$ be objects in $\Rep D(\alpha,\beta)$. They give rise to objects
$(D,\zeta)$ and $(D',\zeta')$ of $Z(\Res^\beta_\alpha)$, by means of the above construction.
On the other hand, since $D(\alpha,\beta)$ is a matched pair,
$D\otimes D'$ has a natural structure of representation of $D(\alpha,\beta)$
(Proposition~\ref{P:monope}), which gives rise to an object $(D\otimes D',\zeta'')$ of
$Z(\Res^\beta_\alpha)$. According to~\eqref{action13}, the transformation 
$\zeta''_{\Ec}:\Res^\beta_\alpha(\Ec)\otimes D\otimes D'\to D\otimes D'\otimes\Res^\beta_\alpha(\Ec)$
is given by
\begin{align*}
\zeta''_{\Ec}(e,d,d') &=
\Bigl(\abs{e}\rightharpoonup(d,d'),\bigl([d,d']\fiz \abs{e}^{-1}\bigr) \fde e\Bigr)\\
&=\Bigl(\abs{e}\rightharpoonup d,\, \bigl(\abs{e}\leftharpoonup(\norm{d},[d])\bigr)\rightharpoonup
d',\, \bigl([d'][d]\fiz\abs{e}^{-1}\bigr)\fde e \Bigr)\\
&=\Bigl(\abs{e}\rightharpoonup d,\, \bigl(([d]\fiz\abs{e}^{-1})\fde\abs{e}\bigr)\rightharpoonup
d',\, \bigl([d'][d]\fiz\abs{e}^{-1}\bigr)\fde e \Bigr)\,.
\end{align*}
In the second equality we made use of~\eqref{action-pt}, \eqref{grading-pt}, and the definition of
the product in $\Vb\prode\Hc^{op}$; in the last equality we used~\eqref{action11bis} (in the special
case when the element of $\Hb$ is an identity).

We need to check that $\zeta''=\zeta\otimes\zeta'$, as prescribed 
by the definition of tensor product of objects in $Z(\Res^\beta_\alpha)$~\eqref{ptcentr}.
We calculate
\begin{align*}
(\zeta\otimes\zeta')_{\Ec}(e,d,d')&\equal{\eqref{ptcentr}}
(\id\otimes\zeta'_{\Ec})(\zeta_{\Ec}\otimes\id)(e,d,d')
\equal{\eqref{action13}}(\id\otimes\zeta'_{\Ec})
\bigl(\abs{e}\rightharpoonup d,\,([d]\fiz\abs{e}^{-1})\fde e,\,d'\bigr)\\
&\equal{\eqref{action13}}
\Bigl(\abs{e}\rightharpoonup d,\,\Abs{([d]\fiz\abs{e}^{-1})\fde e}\rightharpoonup d',\,
\bigl([d']\fiz\Abs{([d]\fiz\abs{e}^{-1})\fde e}^{-1}\bigr)\fde
\bigl(([d]\fiz\abs{e}^{-1})\fde e\bigr)\Bigr)\\
&\equal{\eqref{compcond-action},\eqref{idtres}}
\Bigl(\abs{e}\rightharpoonup d,\,\bigl(([d]\fiz\abs{e}^{-1})\fde\abs{e}\bigr)\rightharpoonup
d',\,
\bigl([d']\fiz([d]\fde\abs{e}^{-1})\bigr)\fde\bigl(([d]\fiz\abs{e}^{-1})\fde e\bigr)\Bigr)\\
&\equal{\eqref{mp-1},\eqref{mp-4}}
\Bigl(\abs{e}\rightharpoonup d,\,\bigl(([d]\fiz\abs{e}^{-1})\fde\abs{e}\bigr)\rightharpoonup
d',\,\bigl([d'][d]\fiz\abs{e}^{-1}\bigr)\fde e\Bigr)
\end{align*}
Thus $\zeta''=\zeta\otimes\zeta'$ and the proof is complete.
\epf

Similarly, one may verify that the rest of the structure (unit object, morphisms) is preserved by
the equivalence. This completes the proof of Theorem~\ref{T:centralizer}.

\begin{obs}\label{R:res-inclusion} According to Theorem~\ref{T:centralizer}, the monoidal
functor $\Rep D(\alpha,\beta)\to\Rep(\Vb,\Hb)$ preserves the forgetful functors to $\Quiv(\Pc)$.
Therefore, by Theorem~\ref{restr-uniq}, this functor must be the restriction along a certain 
morphism of matched pairs $(\Vb,\Hb)\to D(\alpha,\beta)$. It is easy to see that the latter is
the morphism of matched pairs $(\iota,\pi)$ of Remark~\ref{R:inclusion}.
\end{obs}

\begin{corollary}\label{C:dual} Let $(\Vc,\Hc)$ be a matched pair of groupoids with base $\Pc$. There is a
strict monoidal equivalence
\[\Rep(\Hc,\Vc)\cong Z(\ffl)\]
between the category of representations of the dual and the centralizer of the forgetful functor. 
Moreover, the following diagram is commutative
\[\xymatrix@C-20pt@R-5pt{
{\Rep(\Hc,\Vc)}\ar[rr]^{\cong}\ar[rd]_{\ffl}& &{\ \ Z\bigl(\ffl)\ \ \ }\ar[ld]^{U}\\
&{\Quiv(\Pc)}& \\
}\]
\end{corollary}
\pf This is a special case of Theorem~\ref{T:centralizer}, since $\ffl=\Res^{\beta_0}_{\alpha_0}$ where
$(\alpha_0,\beta_0):(\Pc\times\Pc,\Pc)\to(\Vc,\Hc)$ is the unique morphism (Example~\ref{Ex:forgetful}),
and the dual of $(\Vc,\Hc)$ is $D(\alpha_0,\beta_0)$ (Remark~\ref{R:dual}).
\epf

\begin{corollary}\label{C:center} Let $(\Vc,\Hc)$ be a matched pair of groupoids with base $\Pc$.  There is
a strict monoidal equivalence
\[ \Rep D(\Vc,\Hc)\cong Z\bigl(\Rep(\Vc,\Hc)\bigr)\]
between the category of representations of the double and the center of the category of
representations. Moreover, the following diagram is commutative
\[\xymatrix{
{\Rep D(\Vc,\Hc)}\ar[r]^{\cong}\ar[rd]_{\ffl}& {Z\bigl(\Rep(\Vc,\Hc)\bigr)}\ar[r]^{U}&
{\Rep(\Vc,\Hc)}\ar[ld]^{\ffl}\\ 
&{\Quiv(\Pc)}& \\
}\]
\end{corollary}
\pf This is the special case of Theorem~\ref{T:centralizer} when $\Vc = \Vb$, $\Hc = \Hb$, and $\alpha $
and $\beta$ are identities.
\epf

Let $(\Vc,\Hc)$ be a matched pair of groupoids and $D$ a representation of the double 
$D(\Vc,\Hc)=(\Vc\prode\Hc^{op},\Vc\bowtie\Hc)$. Denote the degree of an element $d\in D$ by
$\abs{d}:=(\abs{d}_{\Vc},\abs{d}_{\Hc})\in \Vc\prode\Hc^{op}$ and the action of $(g,h)\in\Vc\bowtie\Hc$
on $d$ by $(g,h)\fde e$. By restriction, this yields actions of $\Vc$ and $\Hc$ on $D$ that we denote
by the same symbol $\fde$. We make use of this notation in the following result. 

\begin{corollary}\label{C:braiding-double} Let $(\Vc,\Hc)$ be a matched pair of groupoids. The category of
representations of the double is braided. If $D$ and $E$ are representations of $D(\Vc,\Hc)$, the braiding
$c_{E, D}: E
\otimes D \to D \otimes E$ is given by
\begin{equation}\label{E:braiding-double}
c_{E, D} (e, d)  =  \bigl(\abs{e}_{\Vc} \fde d, (\abs{d}_{\Hc} \fiz \abs{e}_{\Vc}^{-1}) \fde e\bigr)\,.
\end{equation}
\end{corollary}
\pf This follows from~\eqref{E:braiding-center} and~\eqref{action13}. 
\epf

\section{Braidings on the representations of a matched pair}\label{S:braidings}

We are now in position to classify braidings on the category of representations of a matched
pair in intrinsic terms (that is, in terms of the given groupoid-theoretic structure).

\subsection{Matched pairs of rotations}\label{S:rotations}

\ 

\begin{definition}\label{D:rotation} A {\em rotation} for a matched pair $(\Vc,\Hc)$ is a morphism of
groupoids 
$\kappa:\Vc\to\Hc$ satisfying
\begin{equation}\label{lyz0}
y\kappa(g) = \kappa(y\fde g)(y\fiz g) 
\end{equation}
for all $(g,y)\in \Vc \bowtie \Hc$.
\end{definition}
We may picture a rotation $\kappa$ as follows:
\[\xymatrix@R-10pt@C+5pt{
& {P}\ar@/^/[r]^{y}\ar@/_/[d]_{y\fde g} & {Q}\ar@/^/[d]^{g}\ar@/^/[r]^{\kappa(g)} & {S}\\
{P}\ar@/_/[r]_{\kappa(y\fde g)} &{R}\ar@/_/[r]_{y\fiz g} & {S} & 
}\]

\begin{lema}\label{L:rotation}  Let $(\Vc,\Hc)$ be a matched pair of groupoids and $\kappa: \Vc \to \Hc$ 
a map. Let $\Hat{\kappa}: \Vc \bowtie \Hc \to
\Hc$  be given by $\Hat{\kappa}(g, x) = \kappa(g)x$. Then $\Hat{\kappa}$ is a morphism of groupoids if and 
only if $\kappa$ is a rotation.
\end{lema}
\pf According to~\eqref{E:diag-product}, $\Hat{\kappa}$ is a morphism of groupoids if and only if
$\kappa(f)y\kappa(g)x = \kappa(f(y\fde g))(y\fiz g)x$ for all composable pairs $(f,y)$ and $(g, x)$ in $\Vc
\bowtie\Hc$. This implies that $\kappa: \Vc \to \Hc$ is a morphism of groupoids (letting $y,x$ be
identities) and \eqref{lyz0} (letting $f,x$ be identities). The converse is a special case of the
universal property of the diagonal product~\eqref{E:universal-diagonal}.
\epf

\begin{definition}\label{D:R-matrix} A {\em matched pair of rotations} for a matched pair $(\Vc,\Hc)$
 is a pair of rotations $(\xi, \eta)$  such that 
\begin{align}
\label{lyz0.5} b\bigl(\eta(g) \fde f\bigr) = b\bigl(\xi(f)^{-1} \fde g^{-1}\bigr)\,,\\
\label{lyz1}
\bigl(\eta(g) \fde f\bigr) \bigl(\xi(f)^{-1} \fde g^{-1}\bigr)^{-1}= gf\,.
\end{align}
for every $f$ and $g$ in $\Vc$ with $b(g)=t(f)$.
\end{definition}

The situation may be illustrated as follows:
\[\xymatrix@R-5pt@C+15pt{
& {P}\ar@/_/[d]_{g} & {R}\ar@/_/[l]_{\xi(f)^{-1}\fiz g^{-1}} \\
{P}\ar@/_/[r]^{\eta(g)}\ar@/_/[d]_{\eta(g)\fde f} &{Q}\ar@/_/[d]_{f}\ar@/_/[u]_{g^{-1}} &
{S}\ar@/_/[l]^{\xi(f)^{-1}}\ar@/_/[u]_{\xi(f)^{-1}\fde g^{-1}}\\ 
{R}\ar@/_/[r]_{\eta(g)\fiz f} & {S}\ar@/_/[u]_{f^{-1}} & \\
}\]
This diagram also helps visualize the identities that follow.

\begin{lemma}\label{L:R-matrix} Any matched pair of rotations $(\xi, \eta)$ satisfies the following
conditions:
\begin{align}
\label{lyz2} \xi\bigl(\eta(g) \fde f\bigr) &= \bigl(\xi(f)^{-1} \fiz g^{-1}\bigr)^{-1}\,,\\ 
\label{lyz3}
\eta\bigl(\xi(f)^{-1} \fde g^{-1}\bigr)^{-1} &= \eta(g) \fiz f\,,\\
\label{lyz3.5} \eta(g)\xi(f) &= \bigl(\xi(f)^{-1} \fiz g^{-1}\bigr)^{-1}\bigl(\eta(g) \fiz
f\bigr)\,,
\end{align}
for every $f$ and $g$ in $\Vc$ with $b(g)=t(f)$.
\end{lemma}
\pf Applying the morphism $\xi$ to~\eqref{lyz1} we obtain
\[\xi\bigl(\eta(g) \fde f\bigr) \xi\bigl(\xi(f)^{-1} \fde g^{-1}\bigr)^{-1}= \xi(g)\xi(f)\,.\]
{}From~\eqref{lyz0} for $\xi$ we deduce
\[\xi\bigl(\xi(f)^{-1} \fde g^{-1}\bigr)=\xi(f)^{-1}\xi(g)^{-1}\bigl(\xi(f)^{-1}\fiz g^{-1})^{-1}\,.\]
These two identities combine to give~\eqref{lyz2}.

Similarly, applying the morphism $\eta$ to~\eqref{lyz1} we obtain
\[\eta\bigl(\eta(g) \fde f\bigr) \eta\bigl(\xi(f)^{-1} \fde g^{-1}\bigr)^{-1}= \eta(g)\eta(f)\,.\]
{}From~\eqref{lyz0} for $\eta$ we deduce
\[\eta\bigl(\eta(g) \fde f\bigr) =\eta(g)\eta(f)\bigl(\eta(g)\fiz f\bigr)^{-1}\,.\] 
These two identities combine to give~\eqref{lyz3}.

Finally, from~\eqref{lyz0} for $\xi$ we may also deduce
\[\eta(g)\xi(f)=\xi\bigl( \eta(g)\fde f\bigr)\bigl(\eta(g)\fiz f\bigr)\,.\]
Together with~\eqref{lyz2} this gives~\eqref{lyz3.5}.
\epf

\subsection{Classification of braidings}\label{S:classification}

\  

Braidings on the category of representations of a matched pair are parametrized by matched pairs
of rotations.

\begin{theorem}\label{T:classification}  Let $(\Vc,\Hc)$ be a matched pair of groupoids. There is a
bijective correspondence between matched pairs
of rotations for $(\Vc,\Hc)$ and braidings on $\Rep(\Vc,\Hc)$.
Let $(\xi,\eta)$ be a matched pair
of rotations. For any representations $D$ and $E$, the corresponding braiding 
$c_{E, D}: E \otimes D \to D \otimes E$ is given by
\begin{equation}\label{E:braiding-R-matrix}
c_{E, D} (e, d)  =  \Bigl(\eta(\abs{e}) \fde d, \bigl(\xi(\abs{d})^{-1} 
\fiz \abs{e}^{-1}\bigr) \fde
e\Bigr)\,.
\end{equation}
\end{theorem}

\pf As recalled in Section~\ref{S:generalities}, braidings on $\Rep(\Vc,\Hc)$ are in bijective
correspondence with monoidal sections of the functor $U:Z\bigl(\Rep(\Vc,\Hc)\bigr)\to\Rep(\Vc,\Hc)$.
These are in bijection with sections of 
\[\Rep D(\Vc,\Hc)\map{\cong} Z\bigl(\Rep(\Vc,\Hc)\bigr)\map{U}\Rep(\Vc,\Hc)\,.\]
As noted in Remark~\ref{R:res-inclusion}, this composite is $\Res^\pi_\iota$, where
$\pi: \Vc \prode \Hc^{op} \to \Vc$ and $\iota: \Hc \to \Vc \bowtie \Hc$ are given by 
\[\pi(f, y) = f \text{ \ and \ }\iota(x) = (\id_{\Vc}l(x), x)\,.\]

Note that if a functor $F$ preserves forgetful functors, then so does any section $G$ of $F$:
\[\ffl\circ F=\ffl \quad \Rightarrow \quad \ffl\circ G=\ffl\circ F\circ G=\ffl\,.\]
Therefore, monoidal sections $G$ of $\Res^\pi_\iota$ preserve the forgetful functors and by
Theorem~\ref{restr-uniq} they are in bijective correspondence with morphisms of
matched pairs $(\alpha,\beta): D(\Vc,\Hc)\to (\Vc,\Hc)$ ($\alpha: \Vc \bowtie \Hc \to \Hc$, $\beta: \Vc 
\to \Vc \prode \Hc^{op}$), such that
$(\alpha,\beta)\circ(\iota,\pi)=\id_{(\Vc,\Hc)}$, via $G=\Res^\beta_\alpha$.

By Proposition~\ref{restr} the above condition is equivalent to
\[\alpha\circ\iota = \id_{\Hc} \text{ \ and \ } \pi\circ\beta = \id_{\Vc}\,.\]
In turn, for these conditions to be met, $\alpha$ and $\beta$ must be of the form
\begin{equation*}
\alpha(g, x) = \eta(g)x, \qquad \beta(f) = \bigl(f, \mu(f)\bigr)\,,
\end{equation*}
for certain maps $\mu, \eta: \Vc \to \Hc$ (to obtain the expression for $\alpha$, note that
$(g,x)=\bigl(g,\id_{\Hc} b(g)\bigr)\iota(x)$ in $\Vc\bowtie\Hc$).  

We will show that $(\alpha,\beta)$
is a morphism of matched pairs if and only if $(\xi,\eta)$ is a matched pair
of rotations. This will prove that there
is a bijective correspondence between braidings and matched pairs
of rotations. Expression~\eqref{E:braiding-R-matrix}
for the braiding in $\Rep(\Vc,\Hc)$ follows from expression~\eqref{E:braiding-double} for the
braiding in $\Rep D(\Vc,\Hc)$, taking into account~\eqref{functor3} and~\eqref{functor4} and
the above expressions for $\alpha$ and $\beta$ in terms of $\eta$ and $\xi$.

By Lemma~\ref{L:rotation},
$\alpha$ is a morphism of groupoids if and only if $\eta$ is a rotation. Also,
by~\eqref{E:universal-prode}, $\beta$ is a morphism of groupoids if and only if $\mu:\Vc \to \Hc^{op}$ is
a morphism of groupoids; or equivalently, if $\xi:\Vc\to\Hc$ given by $\xi(f)=\mu(f)^{-1}$ is a morphism
of groupoids. 

Assume then that $\alpha$, $\beta$, $\eta$, and $\xi$ are morphisms of groupoids.
To finish the proof, we will show that conditions~\eqref{functor1} and~\eqref{functor2} for
$(\alpha,\beta)$   are equivalent to condition~\eqref{lyz0} for $\xi$ and conditions~\eqref{lyz0.5}
and~\eqref{lyz1} for $(\xi,\eta)$ (this will complete the proof of the fact that
$(\alpha,\beta)$ is a morphism of matched pairs if and only if $(\xi,\eta)$ is a matched pair
of rotations).
It will turn out that, in this particular case, \eqref{functor2} is a consequence of~\eqref{functor1}.

Let us spell out~\eqref{functor1} and~\eqref{functor2}. They are, respectively,
\[ \beta\bigl(\alpha(g,x)\fde f\bigr) = (g,x)\rightharpoonup \beta(f) \text{ \ and \ }
\alpha\bigl( (g,x)\leftharpoonup(f,\mu(f)\bigr) = \alpha(g,x)\fiz f\,,\]
for every $(g,x)\in\Vc\bowtie\Hc$ and $f\in\Vc$ for which the products are defined.
The actions $\leftharpoonup$ and $\rightharpoonup$ of the matched pair $D(\Vc,\Hc)$ are defined
by~\eqref{action11} and~\eqref{action11bis} (the morphisms $\alpha$ and $\beta$ appearing in those
formulas are in this case identities). Combining these formulas with the expressions for $\alpha$ and
$\beta$ in terms of $\eta$ and $\mu$ we obtain that~\eqref{functor1} is equivalent to
(A) and (B) below, while~\eqref{functor2} is equivalent to (C). 

\begin{flalign}  \tag{A}
\bigl(\eta(g) x\bigr) \fde f &= g(x\fde f)\Bigl(\bigl((x\fiz f) \mu(f)x^{-1}\bigr) \fde g^{-1}\Bigr)\,,\\
\tag{B}
\mu \Bigl(\bigl(\eta(g) x\bigr) \fde f \Bigr)&= \bigl( (x \fiz f)\mu(f) x^{-1} \bigr) \fiz g^{-1}\,,\\
\tag{C}
\bigl(\eta(g) x\bigr) \fiz f &= 
\eta\Biggl(\Bigl(\bigl((x\fiz f)\mu(f)x^{-1}\bigr)\fiz g^{-1}\Bigr)\fde g\Biggr) (x\fiz f) \,,
\end{flalign}
for every $f, g\in \Vc$ and $x\in \Hc$ with $b(g) = l(x)$ and $r(x) = t(f)$.

Letting $x$ be an identity in (A), we get~\eqref{lyz0.5} and~\eqref{lyz1}.
Letting $g$ be an identity in (B), we get that $\xi$ satisfies \eqref{lyz0} (recall $\xi(f)=\mu(f)^{-1}$).
This proves one implication. 

Conversely, assume that $\xi$ satisfies~\eqref{lyz0} and conditions~\eqref{lyz0.5} and~\eqref{lyz1} hold.
By Lemma~\ref{L:R-matrix}, conditions~\eqref{lyz2} and~\eqref{lyz3} hold as well. We then have
\begin{align*} 
\bigl(\eta(g) x\bigr) \fde f &\equal{\eqref{mp-1}} \eta (g)\fde(x\fde f)
\equal{\eqref{lyz1}} g(x\fde f)\bigl(\mu(x\fde f) \fde g^{-1}\bigr)\\ 
&\equal{\eqref{lyz0}} g(x\fde f)\Bigl(\bigl( (x\fiz f)
\mu(f) x^{-1} \bigr)\fde g^{-1} \Bigr)\,,
\end{align*}
which is (A). Also,
 \[\mu \Bigl(\bigl(\eta(g) x\bigr) \fde f \Bigr)\equal{\eqref{mp-1}}
\mu \bigl(\eta (g)\fde(x\fde f) \bigr) \equal{\eqref{lyz2}}
 \mu(x\fde f)\fiz g^{-1} \equal{\eqref{lyz0}}
 \bigl( (x \fiz f)\mu(f) x^{-1} \bigr) \fiz g^{-1}\,,\]
which is (B). Finally,
\begin{align*} 
\bigl(\eta(g) x\bigr) \fiz f  &\equal{\eqref{mp-4}}\bigl(\eta(g)\fiz (x\fde f)\bigr)(x\fiz f)
\equal{\eqref{lyz3}}\eta\Bigl(\bigl(\mu(x \fde f) \fiz g^{-1}\bigr) \fde g\Bigr)(x\fiz f)\\
&\equal{\eqref{lyz0}}\eta\Biggl(\Bigl(\bigl((x\fiz f)\mu(f)x^{-1}\bigr)\fiz g^{-1}\Bigr)\fde g\Biggr)
(x\fiz f)\,,
\end{align*}
which is (C). The proof is complete.
\epf

\begin{exa}\label{Ex:R-matrix-double} Let $(\Vc,\Hc)$ be a matched pair. By
Corollary~\ref{C:braiding-double}, there is a braiding in the category of representations of $D(\Vc,\Hc)$.
According to Theorem~\ref{T:classification}, this braiding must come from a matched pair
of rotations for $D(\Vc,\Hc)$. It follows from~\eqref{E:braiding-double} and~\eqref{E:braiding-R-matrix}
that this consists of the rotations $\xi,\eta:\Vc\prode\Hc^{op}\to\Vc\bowtie\Hc$ given by
\[\xi(\gamma,x)=\bigl(\id_{\Vc}l(x),x^{-1}\bigr) \text{ \ and \ } \eta(\gamma,x)=\bigl(\gamma,\id_{\Hc}
t(\gamma)\bigr)\,.\]
\end{exa}

\section{Applications to weak Hopf algebras}\label{S:applications}

\subsection{Weak Hopf algebras from matched pairs of groupoids}\label{S:weak}

\

We first recall the definition of weak Hopf algebras; more information can be found in~\cite{nik-v}.

\begin{definition}\label{qg}  \cite{bnsz, bsz} A \emph{weak bialgebra} is a collection 
$(H, m, 1, \Delta, \varepsilon)$, where $(H, m, 1)$ is a unital associative algebra, $(H, \Delta, \varepsilon)$
is a counital coassociative coalgebra, and the following conditions hold, for all $a,b,c\in H$:
\begin{gather*}
\Delta(ab) = \Delta(a) \Delta(b)\,, \\ 
\bigl( \Delta(1) \otimes 1 \bigr)\bigl( 1 \otimes \Delta(1) \bigr) = 
\Delta^{(2)} (1 ) = \bigl( 1 \otimes \Delta(1) \bigr)\bigl( \Delta(1) \otimes 1\bigr)\,,\\ 
 \varepsilon(ab_1)\varepsilon(b_2c) =\varepsilon(abc) = \varepsilon(ab_2)\varepsilon(b_1c)\,. 
\end{gather*}
Consider the maps $\varepsilon_{\sou}:H\to H$ and $\varepsilon_{\tgt}(h):H\to H$ defined by
\[\varepsilon_{\sou}(h)=(\id \otimes \varepsilon) \bigl( (1 \otimes h) \Delta(1) \bigr)
\text{ \and \ }
\varepsilon_{\tgt}(h)=(\varepsilon \otimes \id) \bigl(  \Delta(1) (h \otimes 1)\bigr)\,.\]
They are not morphisms of algebras in general, but their images are subalgebras of $H$. 
They are called the source and target subalgebras, respectively. 
Elements of one commute with elements of the other.

A weak bialgebra $H$ is called a \emph{weak Hopf algebra} or a \emph{quantum groupoid}
if there exists a linear map $\Ss: H \to H$ such that
\[ m(\Ss \otimes \id) \Delta = \varepsilon_{\sou}\,, \quad
m(\id \otimes \Ss) \Delta =\varepsilon_{\tgt}\,, \text{ \ and \ }
 m^{(2)}(\Ss \otimes \id \otimes \Ss) \Delta^{(2)}= \Ss\,.\]
\end{definition}
In this case, the map $\Ss$ is unique and it is an antimorphism of algebras
and coalgebras. It is called the antipode of $H$. If $H$ is finite dimensional,
$\Ss$ is invertible and restricts to an anti-isomorphism between the source 
and target subalgebras.

A construction of a weak Hopf algebra from a matched pair of groupoids
was introduced in~\cite{AN}. It generalizes a  construction of Hopf algebras from  matched
pairs of groups due to G. I. Kac and M. Takeuchi. In order to recall it, we need some
terminology from the theory of double categories.

Let $(\Vc,\Hc)$ be a matched pair of groupoids over $\Pc$. The set of cells of $(\Vc,\Hc)$ is
\[\Bg:=\Hc\Times{r}{t}\Vc=\{(x,g) : x\in\Hc,\ g\in\Vc, \text{ and }r(x)=t(g)\}\,.\]
We use the letters $A$, $B$, $C$ to denote cells. A cell may at times be represented by a
diagram~\eqref{E:cell}, which displays the pair $(x,g)$ together with the actions $x\fde g$ and $x\fiz g$.
Two cells $A=(x,g)$ and $B=(y,f)$ are horizontally composable if $g=y\fde f$; in this case, their
horizontal composition is the cell $A B:=(xy,f)$. 
Given an element $g\in\Vc$, we let $I_{\Hc}(g):=\bigl(\id_{\Hc}t(g),g\bigr)$. This is an identity for the
horizontal composition. Vertical composition is defined similarly, and given $x\in\Hc$ we let
$I_{\Vc}(x):=\bigl(x,\id_{\Vc}r(x)\bigr)$. The inverse of a cell $A=(x,g)$ is
$A^{-1}:=\bigl((x\fiz g)^{-1},(x\fde g)^{-1}\bigr)$, as illustrated in~\eqref{E:inversecell}.

We assume from now on that $\Hc$, $\Vc$ (and hence also $\Pc$) are finite sets. The weak Hopf
algebra $\ku(\Vc,\Hc)$ is the $\ku$-vector space with basis $\Bg$ and the following structure,
which we specify on basis elements $A$, $B$:
 
\begin{itemize}
\item The multiplication is 
$A\cdot B = \begin{cases} A B & \text{if $A$ and $B$ are horizontally composable,} \\
0  & \text{otherwise.} \end{cases}$

\medbreak
\item The unit is $\uno := \sum_{g \in \Vc} I_{\Hc}(g)$. 

\medbreak
\item The comultiplication is  $\Delta(A) = \sum B\otimes C$, 
where the sum is over all pairs $(B,C)$ of vertically composable cells whose
vertical composition is $A$.

\medbreak
\item The counit is  $\varepsilon(A) =\begin{cases} 1 & \text{if $A=I_{\Vc}(x)$ for some $x\in\Vc$,}  \\ 
0 & \text{otherwise.}\end{cases}$

\medbreak
\item The antipode is  $\Ss(A)=A^{-1}$.

\end{itemize}

For each $P\in\Pc$, let
\[\uno_{P} := \sum_{g\in\Vc,\, b(g)=P} I_{\Hc}(g) \text{ \ and \ } \uno^{P}:= 
\sum_{g \in \Vc,\,t(g)=P}I_{\Hc}(g)\,.\]
The source subalgebra is the subspace generated by  $\{\uno_{P}\}_{P\in\Pc}$, and
the target subalgebra is the subspace generated by  $\{\uno^{P}\}_{P\in\Pc}$.
Each of these is a complete set of orthogonal idempotents. 
Thus, both subalgebras are semisimple commutative of dimension $\vert \Pc \vert$,
and hence isomorphic to the algebra of functions $\ku^{\Pc}$.

\begin{obs} Two points regarding notation. In the terminology of~\cite{AN}, $\ku(\Vc,\Hc)$  is the weak
Hopf algebra corresponding to the {\em vacant double
groupoid} associated to the transpose of the matched pair $(\Vc,\Hc)$. 
Also, our 
convention for writing products in a groupoid from
left to right is the opposite of that in~\cite{bnsz,nik-v}; for this reason, the source and target
subalgebras of $\ku(\Vc,\Hc)$ are respectively described by the target $b$ and the source $t$ of $\Vc$.

Note that the set $\Bg$ is essentially the same as the set of arrows of the diagonal
groupoid $\Vc\bowtie\Hc$. The algebra $\ku(\Vc,\Hc)$, on the other hand, is {\em not}
the groupoid algebra of $\Vc\bowtie\Hc$, in general.
As an algebra, $\ku(\Vc,\Hc)$ is the groupoid algebra of a groupoid 
whose arrows are the cells of $(\Vc,\Hc)$ and whose objects are the
arrows of $\Vc$. As a coalgebra, it is the dual of the groupoid algebra
of a groupoid whose arrows are the cells of $(\Vc,\Hc)$ and whose 
objects are the arrows of $\Hc$. For more in this direction,
see~\cite{AN}.
\end{obs}

\begin{exa}\label{Ex:weak}
The weak Hopf algebra corresponding to the initial matched pair is the algebra of functions on
$\Pc\times\Pc$:
\[\ku(\Pc\times\Pc,\Pc)\cong \ku^{\Pc\times\Pc}\,.\]
The isomorphism identifies a cell $\bigl(P,(P,Q)\bigr)$ with the function 
$(R,S)\mapsto\begin{cases} 1 & \text{ if $(R,S)=(P,Q)$, }\\
0 & \text{ if not.}\end{cases} $

On the other hand, the weak Hopf algebra corresponding to the terminal matched pair is 
the algebra of linear endomorphisms of the $\ku$-space with basis $\Pc$ (a matrix algebra):
\[\ku(\Pc,\Pc\times\Pc)\cong \End_{\ku}(\ku\Pc)\,.\]
The isomorphism identifies a cell $\bigl((P,Q),Q)\bigr)$ with the endomorphism 
$R\mapsto\begin{cases} P & \text{ if $R=Q$, }\\
0 & \text{ if not.}\end{cases} $

The weak Hopf algebras corresponding to the matched pair $(\Pc,\Hc)$ is the groupoid
algebra of $\Hc$, while the weak Hopf algebra corresponding to $(\Vc,\Pc)$ is the 
algebra of functions $\ku^\Vc$. 
\end{exa}

The construction is functorial with respect to morphisms of matched pairs (Definition~\ref{D:morphism}).

\begin{proposition}\label{P:morph-weak}
Let $(\alpha,\beta):(\Vb,\Hb)\to(\Vc,\Hc)$ be a morphism of matched pairs of finite groupoids.
There is  a morphism of weak Hopf algebras $\ku(\Vb,\Hb)\to\ku(\Vc,\Hc)$ given by
\begin{equation}\label{E:morph-weak}
(h,\gamma)\mapsto
\sum_{g\in\Vc,\,\beta(g)=\gamma}\bigl(\alpha(h),g\bigr)\,.
\end{equation} This defines a functor from the
category of matched pairs to that of weak Hopf algebras.
\end{proposition}
\pf First note that if $(h,\gamma)$ is a cell in $(\Vb,\Hb)$, then $\bigl(\alpha(h),g\bigr)$
is a cell in $(\Vc,\Hc)$ for any $g\in\Vc$ with $\beta(g)=\gamma$, since $\beta:\Vc\to\Vb$ is a morphism
of groupoids. Thus,~\eqref{E:morph-weak} is well-defined.

Let $(h,\gamma)$ and $(h',\gamma')$ be two horizontally composable cells in $(\Vb,\Hb)$. There is a
bijection between the set of cells in $(\Vc,\Hc)$ of the form $\bigl(\alpha(hh'),g'\bigr)$
with $\beta(g')=\gamma'$ and the set of pairs of cells in $(\Vc,\Hc)$ of the form
$\bigl(\alpha(h),g\bigr)$, $\bigl(\alpha(h'),g'\bigr)$, with $\beta(g)=\gamma$ and $\beta(g')=\gamma'$.
These pairs are necessarily horizontally composable. In one direction, the bijection is given by horizontal
composition; in the other, one sets $g:=\alpha(h')\fde g'$ and appeals to~\eqref{functor1}. This shows
that~\eqref{E:morph-weak} preserves multiplications. The remaining conditions can be verified similarly.
\epf
We denote the morphism of matched pairs and the corresponding morphism of weak Hopf algebras by the same
symbol $(\alpha,\beta)$. 

In view of Example~\ref{Ex:weak}, the morphisms from the initial and to the terminal
matched pair (Proposition~\ref{P:initial}) yield morphisms of weak Hopf algebras 
\[\ku^{\Pc\times\Pc}\to \ku(\Vc,\Hc)\to\End_{\ku}(\ku\Pc)\,.\]
The image of the first map is the subalgebra of $\ku(\Vc,\Hc)$ generated by
the source and target subalgebras (which commute elementwise).

\begin{obs}\label{R:basic-weak} The weak Hopf algebra
associated to the opposite of a matched pair $(\Vc,\Hc)$ (Remark~\ref{R:basic}) is the opposite of the
weak Hopf algebra associated to $(\Vc,\Hc)$. Copposites are similarly preserved.
Recall that any matched pair is isomorphic to its opposite and to its coopposite (Remark~\ref{R:iso-opp}).
In view of Proposition~\ref{P:morph-weak}, the same is true of the weak Hopf algebra of a matched pair.
\end{obs}

\subsection{The linearization functor}\label{S:linearization}

 \

Let $H$ be a weak bialgebra. The category $\Repk(H)$ of left $H$-modules 
is a monoidal category with tensor product 
\[M \botimes N := \Delta(1)\cdot M\otimes_{\ku} N\,,\]
where we let $H\otimes H$ act on $M\otimes_{\ku} N$ by $(a,b)\cdot(m,n):=a\cdot m\otimes b\cdot n$.
The action of $H$ on $M \botimes N$ is by restriction of this action via $\Delta$.

Let $(\Vc,\Hc)$ be a matched pair of groupoids. Given a representation $\Ec$ of $(\Vc,\Hc)$,
consider the vector space $\ku\Ec$ with basis $\Ec$. Define a map 
$\ku(\Vc,\Hc) \otimes\ku\Ec \to \ku\Ec$ on the basis $\Bg\times\Ec$ by 
\begin{equation}\label{E:linearaction}
(x,g)\cdot e = \begin{cases} x\fde e & \text{ if $\abs{e}= g$,} \\ 
0 & \text{ otherwise.} \end{cases}
\end{equation}

\begin{proposition}\label{P:linearization} In the  situation above, $\ku\Ec$ is a left module over
$\ku(\Vc,\Hc)$. Moreover, this defines a monoidal functor
\[\Lin:\Rep(\Vc,\Hc)\to\Repk\bigl(\ku(\Vc,\Hc)\bigr)\,.\]
\end{proposition}
\pf Note that if $\abs{e}= g$ and $(x,g)\in\Bg$ then $p(e)=t(g)=r(x)$ so $x\fde e$ is defined.
 Associativity and unitality of~\eqref{E:linearaction} follow from~\eqref{mp-0},~\eqref{mp-1},
and~\eqref{mp-1.5}. A morphism of representations $\psi:\Ec\to\F$ preserves
the actions of $\Hc$ and gradings by $\Vc$, so its linearization
$\psi:\ku\Ec\to\ku\F$ preserves the actions of $\ku(\Vc,\Hc)$. Thus, $\Lin$ is
a functor.

 Let $\Ec$ and $\F$ be representations of $(\Vc,\Hc)$. We identify
 $\ku\Ec\otimes_{\ku}\ku\F$ with $\ku(\Ec\times\F)$. An element 
 $I_{\Hc}(g)\otimes I_{\Hc}(g')$ acts on a basis element $(e,f)$ of 
 $\ku\Ec\otimes_{\ku}\ku\F$ as zero unless $g=\abs{e}$ and $g'=\abs{f}$,
 in which case it acts as the identity, by~\eqref{mp-1.5}. Therefore,
 \[(\uno_P\otimes\uno^P)\cdot(e,f)=
 \begin{cases} (e,f) &\text{ if $b(\abs{e})=P=t(\abs{f})$,}\\
 0 & \text{ otherwise,}\end{cases} \text{ \ and \ }
 \Delta(\uno)\cdot(e,f)=
 \begin{cases} (e,f) &\text{ if $b(\abs{e})=t(\abs{f})$,}\\
 0 & \text{ otherwise,}\end{cases}\]
 since $\Delta(\uno)=\sum_{P\in\Pc}\uno_P\otimes\uno^P$.
Recalling that $\Ec\otimes\F=\{(e,f)\in\Ec\times\F :  b(\abs{e})=t(\abs{f})\}$,
we conclude that
\[\Delta(\uno)\cdot\ku\Ec\otimes_{\ku}\ku\F=\ku(\Ec\otimes\F)\,,\]
or in other words, $\Lin(\Ec)\botimes\Lin(\F)=\Lin(\Ec\otimes\F)$ as vector
spaces. It remains to verify that the actions of $\ku(\Vc,\Hc)$ on both sides
agree. Choose a cell $(x,g)\in\Bg$. 
Consider first its action on a basis element $(e,f)$ of $\ku\Ec\botimes\ku\F$.
For a tensor $B\otimes C$ to appear in
$\Delta(x,g)$ we must have $B=(x,g_1)$ and $C=(x\fiz g_1,g_2)$ with $g_1g_2=g$.
The action of $B\otimes C$ on $(e,f)$ is zero unless $g_1=\abs{e}$ and
$g_2=\abs{f}$, in which case it is equal to 
$\bigl(x\fde e,\,(x\fiz g_1)\fde f\bigr)$. Therefore,
\[\Delta(x,g)\cdot(e,f)=\begin{cases}
\bigl(x\fde e,\,(x\fiz \abs{e})\fde f\bigr) & \text{ if }\abs{e}\abs{f}=g\,,\\
0 & \text{ otherwise.}\end{cases}\]
This agrees with the action of $(x,g)$ on $\ku(\Ec\otimes\F)$,
according to~\eqref{action-pt},~\eqref{grading-pt}, and~\eqref{E:linearaction}.
\epf
\begin{exa}\label{Ex:regularrep}
Consider the distinguished representation $\Hc\Times{r}{t}\Vc$ of a matched
pair $(\Vc,\Hc)$ (Section~\ref{S:distinguished}). It follows
readily from~\eqref{E:dist-rep} and~\eqref{E:linearaction} that the module 
$\ku(\Hc\Times{r}{t}\Vc)$
affords the regular representation of the weak Hopf algebra $\ku(\Vc,\Hc)$.
\end{exa}

The linearization functor satisfies the following naturality condition.
Given a morphism of matched pairs $(\alpha,\beta):(\Vb,\Hb)\to(\Vc,\Hc)$,
the following diagram commutes:
\[\xymatrix@C+10pt{
{\Rep(\Vc,\Hc)}\ar[r]^(0.45){\Lin}\ar[d]_{\Res_\alpha^\beta} & 
{\Rep_{\ku}\bigl(\ku(\Vc,\Hc)\bigr)}\ar[d]^{\Res_\alpha^\beta}\\
{\Rep(\Vb,\Hb)}\ar[r]_(0.45){\Lin} & {\Rep_{\ku}\bigl(\ku(\Vb,\Hb)\bigr)}
}\]
where the functor on the right is the restriction along the morphism of
algebras $\ku(\Vb,\Hb)\to\ku(\Vc,\Hc)$ of Proposition~\ref{P:morph-weak}.

Consider the case of the unique morphism $(\Pc\times\Pc,\Pc)\to(\Vc,\Hc)$.
The category of representations of the initial matched pair is $\Quiv(\Pc)$, 
the category of quivers over $\Pc$ (Example~\ref{Ex:representations}). 
On the other hand, since
$\ku(\Pc\times\Pc,\Pc)=\ku^{\Pc\times\Pc}$ (Example~\ref{Ex:weak}), 
the category of representations of this weak Hopf algebra can be identified with 
$\Bimod(\ku^\Pc)$, the category of bimodules over $\ku^\Pc$. This may also be
described as the category of $(\Pc\times\Pc)$-graded  vector
spaces (with homogeneous maps of degree $0$ as morphisms). We obtain a commutative
diagram
\[\xymatrix@C+10pt{
{\Rep(\Vc,\Hc)}\ar[r]^(0.45){\Lin}\ar[d]_{\ffl} & 
{\Rep_{\ku}\bigl(\ku(\Vc,\Hc)\bigr)}\ar[d]^{\ffl}\\
{\Quiv(\Pc)}\ar[r]_(0.45){\Lin} & {\Bimod(\ku^{\Pc})}
}\]
The functor on the right is the restriction along the
canonical morphism $\ku^{\Pc\times\Pc}\to \ku(\Vc,\Hc)$.
Explicitly, if $M$ is a left $\ku(\Vc,\Hc)$-module, the 
$(\Pc\times\Pc)$-grading on $M$ is defined by
$M_{P,Q} := \uno^Q\cdot(\uno_P\cdot M)$.

\subsection{Quasitriangular weak Hopf algebras}\label{S:quasitriangular}

\

\begin{definition}\label{qtqg}  \cite{bsz,NTV,nik-v}.
A \emph{quasitriangular weak Hopf algebra} is a pair $(H, \R)$
where $H$ is a weak Hopf algebra, $\R\in \Delta^{cop}(1) (H\otimes_k H)\Delta(1)$ 
satisfies
\begin{align}
\label{Runo} \Delta^{cop}(h) \R  &= \R \Delta(h)\,, \text{ \ for all \ } h\in H\,,\\ 
\label{Rdos} (\id \otimes \Delta) (\R) &=\R_{13}\R_{12}\,,\\ 
\label{Rtres} (\Delta \otimes \id) (\R) &=\R_{13}\R_{23},
\end{align}
and there exists an element $\Rs\in \Delta(1) (H\otimes_{\ku} H)\Delta^{cop}(1)$ with
\begin{equation}
\label{Rcuatro}  \R \Rs = \Delta^{cop}(1) \text{ \ and \ } \Rs\R = \Delta(1)\,.
\end{equation}
\end{definition}

In this case, the element $\Rs$ unique, and we say that $\R$ is a quasitriangular
structure on $H$.

\medbreak 

A quasitriangular structure $\R$ on $H$ gives rise to a braiding $c$ in the monoidal
category $\Rep_{\ku}(H)$~\cite[Proposition 5.2.2]{nik-v}. If 
$\R=\sum_i R_i\otimes R'_i$, and $M$ and $N$ are $H$-modules,
the transformation $c_{M,N}:M\botimes N\to N\botimes M$ is given by
\begin{equation}\label{E:basicbraiding}
c_{M,N} (m\otimes n)=\sum_i R'_i\cdot n\otimes R_i\cdot m\,.
\end{equation}
Conversely, suppose that $H$ is a weak Hopf algebra, $(\C,c)$ is a braided category, 
$F:\C\to\Rep_{\ku}(H)$ is a monoidal functor, and there are an object $X$ in $\C$ 
such that $F(X)=H$ and
an element $\R\in H\otimes H$ such that $F(c_{E,D}):F(E)\botimes F(D)\to 
F(D)\botimes F(E)$ is given by~\eqref{E:basicbraiding}, for all objects
$E$ and $D$ of $\C$. Then $\R$ is
a quasitriangular structure on $H$. The proof of~\cite[Proposition XIII.1.4]{kas}
can be easily adapted to derive this result. Briefly,~\eqref{Runo} follows from
the $H$-linearity of $F(c_{X,X})$,~\eqref{Rdos} and~\eqref{Rtres} follow from the
braiding axioms, and~\eqref{Rcuatro} from the invertibility of $F(c_{X,X})$.

We make use of this fact to deduce that matched pairs
of rotations for a matched pair yield 
quasitriangular structures for the corresponding weak Hopf algebra.

\begin{theorem}\label{T:positiveR}  
Let $(\xi,\eta)$ be a matched pair of rotations for $(\Vc,\Hc)$. Define $\R_{\xi, \eta} \in 
\ku(\Vc,\Hc) \otimes_{\ku} \ku(\Vc,\Hc)$  by
\begin{equation}\label{trenzapos}
\R_{\xi,\eta} :=  \sum_{(f,g)\in \Vc\Times{b}{t}\Vc} \bigl(\xi(f)^{-1}\fiz g^{-1},g\bigr)
\otimes \bigl(\eta(g),f\bigr)\,.
\end{equation}
Then $\bigl(\ku(\Vc,\Hc), \R_{\xi,\eta}\bigr)$ is a 
quasitriangular weak Hopf algebra.
\end{theorem}
\pf Let $H=\ku(\Vc,\Hc)$ and $X=\Hc\Times{r}{t}\Vc$. Consider the linearization functor 
$\Lin:\Rep(\Vc,\Hc)\to\Rep_{\ku}(H)$ (Proposition~\ref{P:linearization}).
As mentioned in Example~\ref{Ex:regularrep}, we have $\Lin(X)=H$.
According to the previous remark, it suffices to check that for any
representations $\Ec$ and $\D$ of $(\Vc,\Hc)$ we have
 \[c_{\Ec,\D}(e,d)=\sum_{(f,g)\in \Vc\Times{b}{t}\Vc} 
\bigl(\eta(g),f\bigr)\cdot d\otimes
 \bigl(\xi(f)^{-1}\fiz g^{-1},g\bigr)\cdot e\,.\]
But this is immediate from~\eqref{E:braiding-R-matrix} and 
~\eqref{E:linearaction} (recall that if $(e,d)\in\Ec\otimes\D$ then
$b(\abs{e})=t(\abs{d})$).
\epf

The {\em Drinfeld element} of a quasitriangular weak Hopf algebra $(H,\R)$ is
$u:=\sum_i\Ss(R'_i)R_i$. It satisfies many remarkable properties
~\cite[Proposition 5.2.6]{nik-v}. In addition, when $H$ is involutory 
($\Ss^2=\id$) then $u$ is central and $\Ss(u)=u$.

Let $(\xi,\eta)$ be a matched pair of rotations for a matched pair $(\Vc,\Hc)$, and $\R$ the
quasitriangular structure on $\ku(\Vc,\Hc)$ of Theorem~\ref{T:positiveR}.
The Drinfeld element turns out to be
\[u=\sum_{f\in\Vc}\bigl(\varphi(f),f\bigr)\]
where $\varphi(f):=\xi(f)\eta(f)^{-1}$. Note that $\varphi(f)\in\Hc(P,P)$
where $P=t(f)$. This map $\varphi$ is not a morphism of groupoids but
satisfies the following properties (from which the properties of $u$ may
be deduced):
\[\varphi(x\fde f)=x\varphi(f)x^{-1}\,,\ \ 
\varphi(f)\fde f =f\,, \text{ \ and \ }\varphi(f)\fiz f =\varphi(f^{-1})^{-1}\,,\]
for every $f\in\Vc$ and $x\in\Hc$ with $r(x)=t(f)$. One also finds that
\[u^n=\sum_{f\in\Vc}\bigl(\varphi(f)^n,f\bigr)\]
for every $n\in\Z$.

\subsection{Duals and doubles of weak Hopf algebras of matched pairs}\label{S:double-weak}

\ 

Let $(\Vc,\Hc)$ be a matched pair and $(\Hc,\Vc)$ its dual. 
As noted in Remark~\ref{R:basic}, there is
a bijection between the set of cells of these matched pairs, given by transposition:
\[\Hc\Times{r}{t}\Vc\map{\cong}\Vc\Times{b}{l}\Hc\,,\ \  (x,g)\mapsto (x\fde g,x\fiz g)\,.\]
Define a pairing on $\ku(\Vc,\Hc)\otimes\ku(\Hc,\Vc)$ by
\[\langle (x,g),\,(f,y)\rangle:=\begin{cases}
1 & \text{ if $f=x\fde g$ and $y=x\fiz g$,}\\
0 & \text{ otherwise.}\end{cases}\]

This pairing is non-degenerate and it identifies
 the weak Hopf algebra $\ku(\Hc,\Vc)$ with the dual of the weak Hopf 
 algebra $\ku(\Vc,\Hc)$~\cite[Proposition 3.11]{AN}. 

We briefly review the construction of the {\em Drinfeld double} of a finite dimensional weak Hopf 
algebra $H$~\cite[Section 5.3]{nik-v}. The space $H^*\otimes_{\ku} H$ is an algebra with
multiplication
\begin{equation}\label{E:prod-drinfeld}
(\phi\otimes h)\cdot(\phi'\otimes h'):=\langle\Ss(h_1),\,\phi'_1\rangle\phi'_2\phi\otimes h_2h'
\langle h_3,\,\phi'_3\rangle\,.
\end{equation}
The linear span of the elements
\begin{align}
\label{E:target-double} \phi\otimes zh -\langle z,\epsilon_1\rangle\epsilon_2\phi\otimes h\,, & 
\text{ \ $z$ in the target subalgebra of $H$,}\\
\label{E:source-double} \phi\otimes zh -\langle z,\epsilon_2\rangle\epsilon_1\phi\otimes h\,, & 
\text{ \ $z$ in the source subalgebra of $H$,}
\end{align}
is an ideal for this multiplication. The quotient of $H^*\otimes_{\ku} H$ by this ideal is the
Drinfeld double of $H$. It is a quasitriangular weak Hopf algebra. For more details, see~\cite{nik-v}
or~\cite{NTV}.

\begin{theorem}\label{T:drinfeld} There is an isomorphism of quasitriangular weak Hopf algebras
between the Drinfeld double of $\ku(\Vc,\Hc)$ and the weak Hopf algebra of the
matched pair $(\Vc\prode\Hc^{op},\Vc\bowtie\Hc)$. For basis elements $(f,y)\in\ku(\Hc,\Vc)$
and $(x,g)\in\ku(\Vc,\Hc)$, the isomorphism is given by
\begin{equation}\label{E:iso-double}
(f,y)\otimes (x,g)\mapsto \Bigl((f^{-1},x),\,\bigl(g,(x\fiz g)^{-1}(y^{-1}\fiz f^{-1})x\bigr)\Bigr)\,.
\end{equation}
\end{theorem}
\pf Let $H=\ku(\Vc,\Hc)$. Let $(x,g)$ be a cell in $(\Vc,\Hc)$ (a basis element of $H$) and $(f,y)$ a cell
in $(\Hc,\Vc)$ (a basis element of $H^*$). Conditions~\eqref{E:target-double} and~\eqref{E:source-double}
impose the following relations on
$H^*\otimes_{\ku} H$:
\[(f,y)\otimes(x,g)\equiv 0 \text{ \ whenever $t(f)\neq l(x)$ or $b(x\fde g)\neq l(y^{-1}\fiz f^{-1})$.}\]
Therefore, the quotient of $H^*\otimes_{\ku} H$ by these relations can be identified with the
space with basis elements $(f,y)\otimes(x,g)$, where $t(f)= l(x)$ and $b(x\fde g)=l(y^{-1}\fiz
f^{-1})$. These conditions can be depicted as follows:
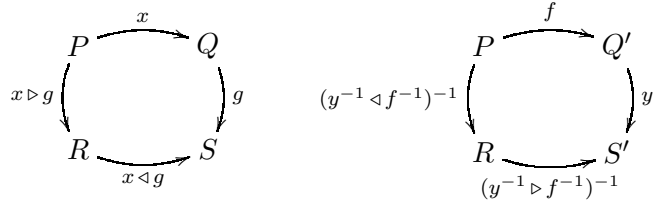
\[\xymatrix@C+10pt{
{P}\ar@/_/[d]_{x\fde g}\ar@/^/[r]^{x} & {Q}\ar@/^/[d]^{g}\\  
{R}\ar@/_/[r]_{x\fiz g} & {S} } \qquad
\xymatrix@C+10pt{
{P}\ar@/_/[d]_{(y^{-1}\fiz f^{-1})^{-1}}\ar@/^/[r]^{f} & {Q'}\ar@/^/[d]^{y}\\  
{R}\ar@/_/[r]_{(y^{-1}\fde f^{-1})^{-1}} & {S'} }\]
We see that $b(f^{-1})=P=l(x)$, so $(f^{-1},x)\in\Vc\bowtie\Hc$, and
$t(g)=Q=r\bigl((x\fiz g)^{-1}(y^{-1}\fiz f^{-1})x\bigr)$ and
$b(g)=S=l\bigl((x\fiz g)^{-1}(y^{-1}\fiz f^{-1})x\bigr)$, so
$\bigl(g,(x\fiz g)^{-1}(y^{-1}\fiz f^{-1})x\bigr)\in\Vc\prode\Hc^{op}$. This shows
that the map~\eqref{E:iso-double} is well-defined. Let us denote it by $\Psi$. Since it takes basis
elements bijectively onto basis elements, $\Psi$ is a linear isomorphism.

Let $\phi=(f,y)$, $h=(x,g)$, $\phi'=(f',y')$, and $h'=(x',g')$.
Computing with~\eqref{E:prod-drinfeld} one finds that the product $(\phi\otimes h)\cdot(\phi'\otimes h')$
is zero unless
\[f'g=(x'\fde g')({y'}^{-1}\fiz{f'}^{-1})^{-1} \text{ \ and \ }
xy'=(y^{-1}\fiz f^{-1})^{-1}(x\fiz g)\,.\]
These conditions are equivalent to
\[\bigl(g,(x\fiz g)^{-1}(y^{-1}\fiz f^{-1})x\bigr)=
({f'}^{-1},x')\rightharpoonup \bigl(g',(x'\fiz g')^{-1}({y'}^{-1}\fiz {f'}^{-1})x'\bigr)\,,\]
as one may see from~\eqref{action11}. This is precisely when the product $\Psi(\phi\otimes
h)\cdot\Psi(\phi'\otimes h')$ is not zero (when these cells are horizontally composable in $D(\Vc,\Hc)$).
Using these conditions one may verify that $\Psi\bigl((\phi\otimes h)\cdot(\phi'\otimes h')\bigr)$
agrees with $\Psi(\phi\otimes h)\cdot\Psi(\phi'\otimes h')$ when these products are non-zero.
We omit the details.
\epf

\section{Appendix. Proof of theorem \ref{T:double}}

We proceed in several steps.
\setcounter{claim}{0}
\begin{claim} The map $\rightharpoonup$  given by~\eqref{action11} is a left action of $\Vc \bowtie \Hb$ on $\Vb \prode \Hc^{op}$.
\end{claim}

\pf Let $(g,h)\in \Vc \bowtie \Hb$, $(\gamma, x) \in \Vb \prode \Hc^{op}$ with $r(h) = t(\gamma)$.
For notational convenience, we set $(A, B)$ for the right-hand side of \eqref{action11}.
Then we also have
$$
A = \beta(g) \left(h\fde \left( \gamma  \beta\left(\left(x\alpha(h)^{-1}\right) \fde g^{-1}\right)\right)
\right).
$$
One may see from~\eqref{E:prod-double} that $A$, $B$, and the action given by~\eqref{action11} are well-defined.
In detail,
 $b(A) = b(\left(\alpha(h\fiz \gamma)x\alpha(h)^{-1}\right) \fde g^{-1}) =   l(B)$ by
\eqref{mp-0.7}; and
$t(A) = t(g) =   r(g^{-1}) = l(B)$ by \eqref{mp-0.3}; hence $(A, B) \in \Vb \prode \Hc^{op}$,
and also \eqref{mp-0} holds.
We check now \eqref{mp-1}. We have
\begin{align*} (f,u) \rightharpoonup ((g,h) \rightharpoonup (\gamma, x)) &=  (f,u) \rightharpoonup  (A, B)
\\ & =
\left(\beta(f) \left(u\fde  \left( A \beta\left(\left( B \alpha(u)^{-1}\right) \fde f^{-1}\right)\right)\right),
\left(\alpha(u \fiz A) B \alpha(u)^{-1}\right) \fiz f^{-1} \right)
\end{align*}
The first component in this expression equals
\begin{equation}\label{aux1}
\beta(f) \left(u\fde  \left( \beta(g) \left( h\fde \left( \gamma  \beta\left(\left(x\alpha(h)^{-1}\right) \fde g^{-1}\right)\right)\right)
\beta\left(\left( B \alpha(u)^{-1}\right) \fde f^{-1}\right)\right)\right),
\end{equation}
while the second is
\begin{equation}\label{aux2}
\left(\alpha(u \fiz \left(\beta(g) \left(h\fde \left( \gamma  \beta\left(\left(x\alpha(h)^{-1}\right) \fde g^{-1}\right)\right) \right)\right))  B \alpha(u)^{-1}\right) \fiz f^{-1}
\end{equation}

On the other hand, the second component of $(f,u)(g,h) \rightharpoonup (\gamma, x) =
(f(\alpha(u)\fde g), (u\fiz \beta(g))h) \rightharpoonup  (\gamma, x)$
equals
\begin{equation}\label{aux3}
\left(\alpha(((u\fiz \beta(g))h)\fiz \gamma)x\alpha((u\fiz \beta(g))h)^{-1}\right) \fiz (f(\alpha(u)\fde g))^{-1},
\end{equation}
while the first is
\begin{align*}
&\quad\;\, \beta(f) \beta((\alpha(u)\fde g) \left((u\fiz \beta(g))h)\fde \left( \gamma
\beta\left(\left(x\alpha(u\fiz \beta(g))h)^{-1}\right) \fde (f(\alpha(u)\fde g))^{-1}\right)\right)
\right),
\\& = \beta(f) (u \fde \beta(g)) \left((u\fiz \beta(g)) \fde \left( h\fde \left( \gamma
\beta\left(\left(x\alpha(u\fiz \beta(g))h)^{-1}\right) \fde (f(\alpha(u)\fde g))^{-1}\right)\right)
\right) \right)
\\& = \beta(f) \left(u \fde \left(\beta(g) \left( h\fde \left( \gamma
\beta\left(\left(x\alpha(u\fiz \beta(g))h)^{-1}\right) \fde (f(\alpha(u)\fde g))^{-1} \right)\right) \right) \right)\right);
\end{align*}
thus, to see that this equals \eqref{aux1}, it is enough to verify the equality
\begin{equation}\label{aux4}
\left( h\fde \left( \gamma  \beta\left(\left(x\alpha(h)^{-1}\right) \fde g^{-1}\right)\right)\right)
\beta\left(\left( B \alpha(u)^{-1}\right) \fde f^{-1}\right)
\overset{?}{=} h\fde \left( \gamma
\beta\left((x\alpha(u\fiz \beta(g))h)^{-1} \fde (f(\alpha(u)\fde g))^{-1} \right)\right).
\end{equation}
The left-hand side of \eqref{aux4} is
\begin{align*}
& \quad \,\; (h\fde \gamma) \left( (h \fiz \gamma ) \fde \left(  \beta\left(\left(x\alpha(h)^{-1}\right) \fde g^{-1}\right)\right)\right)
\beta\left(\left( \left(\left(\alpha(h\fiz \gamma)x\alpha(h)^{-1}\right) \fiz g^{-1} \right) \alpha(u)^{-1}\right) \fde f^{-1}\right)
\\ & = (h\fde \gamma) \beta\left(\left(\alpha(h \fiz \gamma) x\alpha
(h)^{-1}\right) \fde g^{-1}\right) \beta\left(\left( \left(\left(\alpha(h\fiz \gamma)x\alpha(h)^{-1}\right) \fiz g^{-1} \right) \alpha(u)^{-1}\right) \fde f^{-1}\right)
\\ & = (h\fde \gamma) \beta\left(\left(\alpha(h \fiz \gamma) x\alpha
(h)^{-1}\right) \fde \left( g^{-1} \left( \alpha(u)^{-1} \fde f^{-1}\right)\right)\right)
\\&= (h\fde \gamma) \beta\left( \alpha(h \fiz \gamma)\fde \left(
\left( x\alpha (h)^{-1}\right) \fde \left( g^{-1} \left( \alpha(u)^{-1} \fde f^{-1}\right)\right)\right)\right)
\\&= h\fde \left( \gamma \beta\left(
\left( x\alpha (h)^{-1}\right) \fde \left( g^{-1} 
\left( \alpha(u)^{-1} \fde f^{-1}\right)\right)\right) \right)\,.
\end{align*}
Thus, to verify~\eqref{aux4}, it is enough to check
\begin{equation}\label{aux5}
\alpha (h)^{-1} \fde \left( g^{-1} \left( \alpha(u)^{-1} \fde f^{-1}\right)\right)
\overset{?}{=} (\alpha(u\fiz \beta(g))h)^{-1} \fde \left((\alpha(u)\fde g)^{-1} f^{-1} \right).
\end{equation}
Now, the right-hand side of~\eqref{aux5} is
\begin{align*}
& \quad \,\; \left(\alpha(h)^{-1}\alpha(u\fiz \beta(g))^{-1}\right) \fde \left((\alpha(u)\fde g)^{-1} f^{-1} \right)
= \alpha(h)^{-1} \fde \left( (\alpha(u)\fiz g)^{-1}
\fde \left((\alpha(u)\fde g)^{-1} f^{-1} \right)
\right)
\end{align*}
and this equals the left-hand side of \eqref{aux5} by \eqref{idsiete}.
Thus, we have checked the equality of the first components in \eqref{mp-1}. 
We next verify the equality of the second components, that is, the 
equality of \eqref{aux2} and \eqref{aux3}. Acting on \eqref{aux2} by $\fiz f$ 
we have
\begin{align*}
& \quad \,\;
\alpha \left(u \fiz \left(\beta(g) \left(h\fde \left( \gamma  \beta\left(\left(x\alpha(h)^{-1}\right) \fde g^{-1}\right)\right) \right)\right)\right)  \left(\left(\alpha(h\fiz \gamma)x\alpha(h)^{-1}\right) \fiz g^{-1} \right) \alpha(u)^{-1}
\\&= \alpha \left(u \fiz \left(\beta(g) (h\fde \gamma)
\left((h\fiz \gamma) \fde   \beta\left(\left(x\alpha(h)^{-1}\right) \fde g^{-1}\right)\right) \right) \right)  \left(\left(\alpha(h\fiz \gamma)x\alpha(h)^{-1}\right) \fiz g^{-1} \right) \alpha(u)^{-1}
\\&= \alpha \left(u \fiz \left(\beta(g) (h\fde \gamma)
\beta\left(\left(\alpha((h\fiz \gamma))x\alpha(h)^{-1}\right) \fde g^{-1}\right)\right) \right) \left(\left(\alpha(h\fiz \gamma)x\alpha(h)^{-1}\right) \fiz g^{-1} \right) \alpha(u)^{-1}
\\&= \left( \left(\alpha \left(u \fiz \left(\beta(g) (h\fde \gamma) \right) \right) \fiz
\left(\alpha(h\fiz \gamma)x\alpha(h)^{-1}\right) \fde g^{-1}\right)  \right) \left(\left(\alpha(h\fiz \gamma)x\alpha(h)^{-1}\right) \fiz g^{-1} \right) \alpha(u)^{-1}
\\&= \left(\left(\alpha \left(u \fiz \left(\beta(g) (h\fde \gamma) \right)\right)
\alpha(h\fiz \gamma)x\alpha(h)^{-1} \right) \fiz g^{-1}  \right)  \alpha(u)^{-1}
\\&= \left(\left(\alpha \left(u \fiz \left(\beta(g) (h\fde \gamma) \right)
(h\fiz \gamma)\right)x\alpha(h)^{-1} \right) \fiz g^{-1}  \right)  \alpha(u)^{-1}
\\&= \left(\left(\alpha \left(\left(u \fiz \left(\beta(g)(h\fde \gamma) \right) \right)
(h\fiz \gamma)\right)x\alpha(h)^{-1} \right) \fiz g^{-1}  \right)  \alpha(u)^{-1}
\\&= \left(\left( \left(\alpha \left(\left(u \fiz \beta(g) \right)h\right)\fiz \gamma\right)x\alpha(h)^{-1} \right) \fiz g^{-1}\right)  \alpha(u)^{-1}
\end{align*}
On the other hand, \eqref{aux3} acted by $\fiz f$ equals
\begin{align*}
& \quad \,\; \left(\alpha(((u\fiz \beta(g))h)\fiz \gamma)x\alpha((u\fiz \beta(g))h)^{-1}\right) \fiz (\alpha(u)\fde g)^{-1}
\\&= \left(\alpha(((u\fiz \beta(g))h)\fiz \gamma)x\alpha(h)^{-1}
(\alpha(u)\fiz g)^{-1}\right) \fiz (\alpha(u)\fde g)^{-1}
\\&= \left(\left( \alpha(((u\fiz \beta(g))h)\fiz \gamma)x\alpha(h)^{-1}\right) \fiz
\left( (\alpha(u)\fiz g)^{-1} \fde (\alpha(u)\fde g)^{-1}\right)\right)
\left( (\alpha(u)\fiz g)^{-1} \fiz (\alpha(u)\fde g)^{-1}\right)
\\&= \left(\left( \alpha(((u\fiz \beta(g))h)\fiz \gamma)x\alpha(h)^{-1}\right) \fiz g^{-1}\right) \alpha(u)^{-1}.
\end{align*}
Therefore \eqref{mp-1} holds.
The verification of \eqref{mp-1.5} is straightforward, and the claim is proved.
\epf

\begin{claim}
The map $\leftharpoonup$  given by  \eqref{action11bis} is a right action of $\Vb \prode \Hc^{op}$ on $\Vc \bowtie \Hb$. 
\end{claim}

\pf Let $(g,h)\in \Vc \bowtie \Hb$, $(\gamma, x) \in \Vb \prode \Hc^{op}$. Action \eqref{action11bis} is
well-defined since
$$
b\left( ((\alpha(h\fiz \gamma)x\alpha(h)^{-1}) \fiz g^{-1}) \fde g\right)
= l\left( ((\alpha(h\fiz \gamma)x\alpha(h)^{-1}) \fiz g^{-1}) \fiz g\right)= l\left(\alpha(h\fiz \gamma)x\alpha(h)^{-1}\right)= l\left(h\fiz \gamma\right).
$$
We check \eqref{mp-0.3}:
$$
r\left((g,h) \leftharpoonup (\gamma, x)\right) =  r\left( ((\alpha(h\fiz \gamma)x\alpha(h)^{-1}) \fiz g^{-1}) \fde g, h\fiz \gamma\right)
= r(h\fiz \gamma) = b(\gamma) = l(x) = b(\gamma, x).$$
We now check~\eqref{mp-2}. We have $((g,h) \leftharpoonup(\gamma, x) ) \leftharpoonup  (\tau, y) =$
\begin{align*}  &=
\left( ((\alpha(h\fiz \gamma)\, x \, \alpha(h)^{-1}) \fiz g^{-1}) \fde g, h\fiz \gamma\right)
\leftharpoonup (\tau, y)
\\ & = \left(\left((\alpha(h\fiz \gamma) \fiz \tau) \, y \, \alpha(h\fiz \gamma)^{-1})
\fiz ((\alpha(h\fiz \gamma)\, x \,\alpha(h)^{-1}) \fde g^{-1}) \right)
\fde \left( ((\alpha(h\fiz \gamma)\, x \,\alpha(h)^{-1}) \fiz g^{-1}) \fde g \right), \right.
\\ &  \qquad\qquad\qquad\qquad\qquad\qquad\qquad\qquad
\qquad\qquad\qquad\qquad\qquad\qquad\qquad\qquad
\left. (h\fiz \gamma) \fiz \tau \right)
\\ & = \left( \left(\alpha(h\fiz \gamma\tau) \, y \, \alpha(h\fiz \gamma)^{-1}) \fiz ((\alpha(h\fiz \gamma)\, x \,\alpha(h)^{-1}) \fde g^{-1}) ((\alpha(h\fiz \gamma)\, x \,\alpha(h)^{-1}) \fiz g^{-1}\right) \fde g, h\fiz \gamma\tau \right)
\\ & = \left( \left((\alpha(h\fiz \gamma\tau) \, y \, \alpha(h\fiz \gamma)^{-1}\alpha(h\fiz \gamma)\, x \, \alpha(h)^{-1})  \fiz g^{-1}\right) \fde g, h\fiz \gamma\tau \right)
\\ & = \left( \left((\alpha(h\fiz \gamma\tau) \, yx \,\alpha(h)^{-1})  \fiz g^{-1}\right) \fde g, h\fiz \gamma\tau \right)
\\ &= (g,h) \leftharpoonup(\gamma\tau, yx).
\end{align*}
Here we have used $((\alpha(h\fiz \gamma)x\alpha(h)^{-1}) \fiz g^{-1}) 
\fde g)^{-1} = (\alpha(h\fiz \gamma)x\alpha(h)^{-1}) \fde g^{-1}$ in 
the second equality, which follows from \eqref{mp-3}.
Hence, \eqref{mp-2} is valid for \eqref{action11bis}. A straightforward 
computation, using \eqref{iddos}, shows that  also \eqref{mp-2.5} holds.   
\epf

\begin{claim} The groupoids $\Vc \bowtie \Hb$ and $\Vb \prode \Hc^{op}$,
 together with the actions \eqref{action11} and \eqref{action11bis}, 
 form a matched pair. 
\end{claim}

\pf We check \eqref{mp-0.7}:
\begin{align*}
b\left((g,h) \rightharpoonup (\gamma, x)\right) &=
b\left(\beta(g) (h\fde \gamma)  \beta\left(\left(\alpha(h\fiz \gamma)x\alpha(h)^{-1}\right) \fde g^{-1}\right)\right)
\\ & = b\left(\left(\alpha(h\fiz \gamma)x\alpha(h)^{-1}\right) \fde g^{-1}\right)
\\ & = l\left(\left(\alpha(h\fiz \gamma)x\alpha(h)^{-1}\right) \fiz g^{-1}\right);
\\  &=
t\left( ((\alpha(h\fiz \gamma)x\alpha(h)^{-1}) \fiz g^{-1}) \fde g\right)
\\ &= l \left((g,h) \leftharpoonup (\gamma, x)\right).
\end{align*}
We check \eqref{mp-3}:
\begin{align*}
(g,h) \rightharpoonup\left((\gamma, x) (\tau, y)\right) &=  (g,h) \rightharpoonup (\gamma\tau, yx)
\\ & = \left(\beta(g) (h\fde \gamma\tau)  \beta\left(\left(\alpha(h\fiz \gamma\tau)yx\alpha(h)^{-1}\right) \fde g^{-1}\right),
\left(\alpha(h\fiz \gamma\tau)yx\alpha(h)^{-1}\right) \fiz g^{-1} \right),
\end{align*}
while $\left((g,h) \rightharpoonup(\gamma, x)\right)
\left(\left((g,h) \leftharpoonup(\gamma, x)\right)\rightharpoonup (\tau, y)\right)$ equals
\begin{align*}
& \quad\,\;
\left(\beta(g) (h\fde \gamma)  \beta\left(\left(\alpha(h\fiz \gamma)x\alpha(h)^{-1}\right) \fde g^{-1}\right),
\left(\alpha(h\fiz \gamma)x\alpha(h)^{-1}\right) \fiz g^{-1} \right)
\\ &\quad \times
\left(\left( ((\alpha(h\fiz \gamma)x\alpha(h)^{-1}) \fiz g^{-1}) \fde g, h\fiz \gamma\right)\rightharpoonup (\tau, y)\right)
\end{align*}
The first component of this expression is
\begin{align*}
& \quad\,\; \beta(g) (h\fde \gamma)  \beta\left(\left(\alpha(h\fiz \gamma)x\alpha(h)^{-1}\right) \fde g^{-1}\right)
\beta\left(((\alpha(h\fiz \gamma)x\alpha(h)^{-1}) \fiz g^{-1}) \fde g\right) \left((h\fiz \gamma)\fde \tau\right)
\\ &\quad \times \beta\left(\left(\alpha(h\fiz \gamma\tau)y\alpha(h\fiz \gamma)^{-1}\right) \fde \left( (\alpha(h\fiz \gamma)x\alpha(h)^{-1}) \fde g^{-1}\right)\right)
\\ &=  \beta(g) (h\fde \gamma\tau) \beta\left((\alpha(h\fiz \gamma\tau)yx\alpha(h)^{-1}) \fde g^{-1}\right);
\end{align*}
and the second component is
\begin{align*}
& \quad\,\; \left(\left(\alpha(h\fiz \gamma\tau)y\alpha(h\fiz \gamma)^{-1}\right) \fiz \left( (\alpha(h\fiz \gamma)x\alpha(h)^{-1}) \fde g^{-1}\right) \right)
\left(\left(\alpha(h\fiz \gamma)x\alpha(h)^{-1}\right) \fiz g^{-1} \right)
\\ & = \left(\alpha(h\fiz \gamma\tau)y\alpha(h\fiz \gamma)^{-1}\alpha(h\fiz \gamma)x\alpha(h)^{-1}\right) \fiz g^{-1}.
\end{align*}
Hence \eqref{mp-3} holds.
We finally check \eqref{mp-4}. It is convenient to set $Q := \alpha\left(\left((u\fiz \beta(g))h\right) \fiz \gamma\right)x\alpha ( h^{-1})$.  Then the first component of $\left((f,u)(g,h)\right) \leftharpoonup (\gamma, x) = (f(\alpha(u)\fde g), (u\fiz \beta(g))h) \leftharpoonup (\gamma, x)$ is
\begin{align*}
&\quad\,\; \left(\left(\alpha\left(\left((u\fiz \beta(g))h\right) \fiz \gamma)x\alpha\left((u\fiz \beta(g)\right)h\right)^{-1}\right) \fiz \left(f(\alpha(u)\fde g)\right)^{-1}\right) \fde \left(f(\alpha(u)\fde g)\right)
\\ &= \left(\left(\left(Q\left(\alpha(u)\fiz g\right)^{-1}\right) \fiz \left(\alpha(u)\fde g\right)^{-1}\right) \fiz f^{-1}\right)
\fde \left(f(\alpha(u)\fde g)\right)
\\ &= \left(\left(\left(Q\fiz g^{-1}\right) \alpha(u)^{-1}\right) \fiz f^{-1}\right)
\fde \left(f(\alpha(u)\fde g)\right)
\\ &= \left(\left(\left(\left(Q\fiz g^{-1}\right) \alpha(u)^{-1}\right) \fiz f^{-1}\right)\fde f\right)
\left(\left(\left(\left(\left(Q\fiz g^{-1}\right) \alpha(u)^{-1}\right) \fiz f^{-1}\right)\fiz f\right) \fde (\alpha(u)\fde g)\right)
\\ &= \left(\left(\left(\left(Q\fiz g^{-1}\right) \alpha(u)^{-1}\right) \fiz f^{-1}\right)\fde f\right)
\left(\left(\left(Q\fiz g^{-1}\right) \alpha(u)^{-1}\right) \fde (\alpha(u)\fde g)\right)
\\ &= \left(\left(\left(\left(Q\fiz g^{-1}\right) \alpha(u)^{-1}\right) \fiz f^{-1}\right)\fde f\right)
\left(\left(Q\fiz g^{-1}\right) \fde g\right),
\end{align*}
where we have used \eqref{idocho}; while the second is
$\left((u\fiz \beta(g))h\right)\fiz \gamma$.
Set $P= \alpha(h\fiz \gamma)x\alpha(h)^{-1}$. Then
the first component of $((f,u) \fiz ((g,h) \fde (\gamma, x))) ((g,h) \fiz (\gamma, x))$ is
\begin{align*}
&\quad\,\; \left(\left(\alpha \left(u \fiz \left(\beta(g) (h\fde \gamma) \beta(P \fde g^{-1} \right) \right)
\left( (P \fiz g^{-1})\alpha(u)^{-1}\right) \fiz f^{-1} \right) \fde f\right)
\\ &\times \left(\alpha \left(u \fiz \left(\beta(g) (h\fde \gamma) \beta(P \fde g^{-1})\right)  \fde \left((P \fiz g^{-1}) \fde g\right)\right)\right)
\end{align*}
The first factor in the last expression is
\begin{align*}
&\quad\,\; \left(\left(\alpha \left(u \fiz \left(\beta(g) (h\fde \gamma)\right) \fiz (P \fde g^{-1}  \right)\right)
\left((P \fiz g^{-1})\alpha(u)^{-1}\right) \fiz f^{-1} \right) \fde f
\\ &= \left(\left(\left(\left(\alpha \left(u \fiz \left(\beta(g) (h\fde \gamma)\right)\right)P\right)\fiz g^{-1} \right)\alpha(u)^{-1}\right) \fiz f^{-1} \right) \fde f
\\ &= \left(\left(\left(\left(\alpha \left(u \fiz \left(\beta(g) (h\fde \gamma)\right)(h\fiz \gamma)\right)x\alpha(h)^{-1}\right)\fiz g^{-1} \right)\alpha(u)^{-1}\right) \fiz f^{-1} \right) \fde f
\\ &= \left(\left(\left(\left(Q\fiz g^{-1}\right) \alpha(u)^{-1}\right) \fiz f^{-1}\right)\fde f\right);
\end{align*}
while the second factor is
\begin{align*}
&\quad\,\; \left( \alpha \left(u \fiz \left(\beta(g) (h\fde \gamma) \beta(P \fde g^{-1})\right)\right) (P \fiz g^{-1}) \right)\fde g
\\ &= \left( \alpha \left(u \fiz \left(\beta(g) (h\fde \gamma) \right) P\right) \fiz g^{-1} \right)\fde g
\\ &= \left(\left(Q\fiz g^{-1}\right) \fde g\right).
\end{align*}
This shows that the first components of both sides of \eqref{mp-4} in the present setting are equal. Finally,
the second component of $((f,u) \fiz ((g,h) \fde (\gamma, x))) ((g,h) \fiz (\gamma, x))$ is
\begin{align*}
&\quad\,\; \left(\left(u \fiz \left(\beta(g) (h\fde \gamma) \beta(P \fde g^{-1})\right)\right) \fiz \beta\left( (P \fiz g^{-1}) \fde g \right)\right)(h \fiz \gamma)
\\ &= \left(\left(u \fiz \left(\beta(g) (h\fde \gamma) \right)\right) \fiz \beta\left((P \fde g^{-1}) ((P \fiz g^{-1}) \fde g) \right)\right)(h \fiz \gamma)
\\ &= \left(u \fiz \left(\beta(g) (h\fde \gamma) \right)\right) (h \fiz \gamma)
\\ &= \left((u\fiz \beta(g))h\right)\fiz \gamma
\end{align*}
Hence, also the second components of both sides of \eqref{mp-4} in the present setting are equal.
 \epf

The proof of Theorem~\ref{T:double} is complete.

\end{document}